\documentclass{article}

\usepackage[T1]{fontenc}
\usepackage[utf8]{inputenc}
\usepackage{lmodern} %
\usepackage[english]{babel}
\usepackage{hyperref}

\usepackage{xcolor}
\usepackage{graphicx}
\graphicspath{{figures/},{svg-inkscape/},{tikzs/}}
\usepackage{subcaption}

\usepackage{siunitx}
\usepackage{upgreek}

\usepackage{bm}
\usepackage{cancel}
\usepackage{mathtools} %
\usepackage{cleveref}

\let\originalleft\left
\let\originalright\right
\renewcommand{\left}{\mathopen{}\mathclose\bgroup\originalleft}
\renewcommand{\right}{\aftergroup\egroup\originalright}

\newlength\figureheight 
\newlength\figurewidth
\newlength\svgwidth
\newcommand{\includetikz}[1]{\includegraphics{tikzs/#1_compiled}}
\newcommand{\includelegend}[1]{\includegraphics{tikzs/#1}}

\newcommand{\R}{{\mathbb{R}}} 
\newcommand{\tunderbrace}[2]{\underbrace{#1}_{\textstyle#2}}

\newcommand{\hastobe}{\overset{!}{=}}
\DeclareMathOperator{\Div}{div}
\DeclareMathOperator{\trace}{tr}
\DeclareMathOperator{\sym}{sym}

\renewcommand{\d}{\,\mathrm{d}}

\newcommand{\state}{\bm{w}}
\newcommand{\adjstate}{\bm{z}}
\newcommand{\shape}{u}
\newcommand{\fracture}{\mathcal{F}}
\newcommand{\stress}{\bm{\sigma}}
\newcommand{\strain}{\bm{\varepsilon}}
\newcommand{\hmetric}{\mathcal{H}}
\newcommand{\vV}{\bm{V}}
\newcommand{\vW}{\bm{W}}
\newcommand{\deriv}{\mathrm{d}}
\newcommand{\shapeDeriv}{\partial_\remainderdomain}

\newcommand{\physdomain}{\Omega}
\newcommand{\remainderdomain}{U}
\newcommand{\holdalldomain}{\Xi}

\newcommand{\Gb}{\varGamma_\text{bottom}}
\newcommand{\Gl}{\varGamma_\text{left}}
\newcommand{\Gr}{\varGamma_\text{right}}
\newcommand{\Gt}{\varGamma_\text{top}}

\DeclareMathOperator*{\argmin}{arg\,min}

\DeclareMathOperator{\Diff}{Diff}
\newcommand{\Diffc}{\Diff_c}

\newcommand{\Diffcn}{\Diffc(\R^2)}

\newcommand{\smooth}{C^{\infty}}

\usepackage{authblk}
\usepackage{amsmath,amssymb}

\providecommand{\keywords}[1]
{
	\small	
	\textit{Keywords: } #1
}

\begin{document}

\title{Fracture propagation by using
	 	shape optimization techniques based on outer Riemannian metrics}
  
\author[1]{Tim Suchan}
\author[2]{Winnifried Wollner}
\author[3]{Kathrin Welker}

\affil[1]{Helmut Schmidt University/University of the Federal Armed Forces Hamburg, Germany, \texttt{suchan@hsu-hh.de}}
             
\affil[2]{University of Hamburg, \texttt{winnifried.wollner@uni-hamburg.de}}

\affil[3]{Helmut Schmidt University/University of the Federal Armed Forces Hamburg, Germany, \texttt{welker@hsu-hh.de}}

\date{}

\maketitle              %

\begin{abstract}
In this work, we investigate a novel approach for the simulation of two-dimensional,
brittle, quasi-static fracture problems based on a shape optimization approach.
In contrast to the commonly-used phase-field approach, this proposed approach for
investigating fracture paths does not require a specified `length-scale' parameter
defining the diffuse interface region nor a level set function. Instead, it interprets the fracture as part
of the boundary of the domain and uses shape optimization algorithms to minimize the
energy in the system and therefore describes the fracture propagation directly.
Embedding the problem of energy minimization in a Riemannian manifold framework formulated
on a suitable shape space, together with the choice of an outer Riemannian metric, yields
both advantages from an analytical as well as an applied perspective. Furthermore, an
eigenvalue decomposition of the strain tensor is used to produce more realistic fracture
paths (the so-called strain splitting), which only allows fracture growth from tensile loads.
Numerical simulations for the commonly considered single-edge notch tension and shear test
are performed and the results are evaluated in comparison to phase-field results.

\end{abstract}


\section{Introduction}\label{sec:introduction}

Fracture propagation in brittle materials goes back to the model by~\cite{Griffith1921}.
In~\cite{Francfort1998}, a variational formulation for fracture propagation was proposed based
on an energy minimization property for the total energy consisting of a contribution from the
energy due to, elastic, deformation of the unbroken material and a fracture energy,
see also~\cite{Bourdin2008}. To avoid difficulties with the discretization of the lower
dimensional fracture~\cite{BourdinFrancfortMarigo:2000} introduced an approximation by elliptic
functionals following general ideas of~\cite{Ambrosio1990} giving rise to the nowadays common
approximation of fracture propagation by phase-field methods, see,
e.g.,~\cite{Miehe2010a,Ambati2015,Mang2019}. Different models for
tension-compression splitting are considered in literature, see,
e.g.~\cite{Ambati2015,Miehe2010}, to prevent fracture propagation under compression.
However, the elliptic approximation parameter requires a sufficiently finely resolved mesh
and the coupling of these parameters for reliable simulation results requires care, see,
e.g.,~\cite{Wheeler2014}, and specifically designed adaptive discretization schemes,
e.g.,~\cite{BurkeOrtnerSuli:2010,WallothWollner:2021}. To avoid the delicate choice of these
parameters and the interpretation of intermediate phase-field values, within this paper,
we consider a sharp fracture whose evolution is computed by shape optimization techniques.

Shape optimization is a classical topic in mathematics which is of
high importance in a wide range of applications.
In particular, shape optimization problems arise frequently in technological processes modeled
by partial differential equations (PDE), e.g., aerodynamic shape
optimization~\cite{AIAA-2013}, acoustic shape
optimization~\cite{Berggren-horn-2007} and optimization of interfaces
in transmission problems~\cite{Langer-2015}. 
Application-oriented questions in shape optimization are concerned
with specific calculations of shapes which are optimal with respect to
an objective functional dependent on the solution of a PDE.

Furthermore, the topic of fracture propagation has been considered in \cite{Leugering2015}, where shape and topology optimization were used to control the propagation of fracture. There, fracture propagation is modeled using the variational formulation from~\cite{Francfort1998}.
In~\cite{Allaire2011}, the evolution of the fracture itself is considered when a loading is applied. The shape is represented by a level set function, the minimization problem is also based on the variational formulation from~\cite{Francfort1998}, and a shape and topology optimization approach are combined to evolve the level set. A similar approach is followed in~\cite{Alidoost2022}, where also a level set approach is used to determine the evolution of the fracture.

When considering shapes as an open and bounded subset of $\R^2$ or $\R^3$, which is a conceptually different approach compared to using level set functions or phase fields and is followed in this manuscript, then solving shape optimization problems is made more difficult by the fact
that the set of permissible shapes generally does not allow a vector space structure. This concept allows for a less restricted description of the shape, however, in contrast to the previously-mentioned methods, it does not permit the merging of fractures.
In general, the modeling of a shape space and an associated distance
measure for the shapes is a challenging task and different approaches lead to diverse models.
Shapes (and their similarities) have been extensively studied in recent
decades. David G.~Kendall~\cite{Kendall} has already introduced the
notion of a shape space in~1984.
However, there is a large number of different shape concepts, e.g.,
landmark vectors~\cite{Kendall,HafnerZachSand,SoehnBirknerYanAlber},
plane curves~\cite{Michor2006b,Michor2007,MioSrivastavaJoshi},
surfaces~\cite{BauerHarmsMichor,Michor2005}, boundary contours of
objects~\cite{FuchsJuettlerScherzerYang,LingJacobs,RumpfWirth2},
multiphase objects~\cite{WirthRumpf}, characteristic functions of
measurable sets~\cite{Zolesio} and morphologies of images~\cite{DroskeRumpf}.
In an industrial context,  often a priori parametrizations
of the shapes of interest are still used because of the resulting vector
space framework matching standard optimization software. However,
only shapes corresponding to the a priori parametrization can be reached. 
If one cannot work in vector spaces, Riemannian manifolds are the next best option. 
The space
$B_e(S^1,\mathbb{R}^2)\coloneqq \textup{Emb}(S^1,\mathbb{R}^2)/\textup{Diff}(S^1)$
of smooth shapes (cf.~\cite{Michor2006b})  is an important example for a Riemannian manifold.
Here, $S^1$ denotes the one-dimensional unit sphere,
$\textup{Emb}(S^1, \mathbb{R}^2)$ denotes the set of all embeddings of
$S^1$ into $\R^2$ and $\textup{Diff}(S^1)$ denotes the set of all diffeomorphisms of $S^1$. 
This shape space is considered in recent publications (see,
e.g.,~\cite{Rozan,Geiersbach2021,Geiersbach2024,Schulz2014}),
with particular emphasis on an inner Riemannian metric, i.e., a Riemannian metric induced by
right-invariant Riemannian metrics from $\textup{Emb}(S^1, \mathbb{R}^2)$.
Unfortunately, using inner Riemannian metrics, we obtain (shape) gradients which
only deform the shape itself and not the ambient space, i.e., $\R^2$. 
Inspired by the analysis proposed in~\cite{Bauer2011new}, where outer
and inner Riemannian metrics were compared in surface registration
problems,~\cite{LoayzaRomero2025}
proposes to consider outer Riemannian metrics in the context of PDE-constrained
shape optimization instead. An outer Riemannian metric on $B_e$ is given by a
right-invariant Riemannian metric on the space 
\begin{equation}
  \label{Diffcn}
  \Diffcn \coloneqq \{\varphi\in \smooth(\R^2, \R^2) \colon \varphi^{-1} \in \smooth(\R^2, \R^2), \,\operatorname{supp}(\varphi - \text{id}) \text{ is compact}\}
\end{equation}
inducing a Riemannian metric on $\text{Emb}(S^1, \mathbb{R}^2)$ by left action
(cf.~\cite{Bauer2013overview}). Hereby, $\text{id}$ denotes the
identity map, $\text{supp}$ the support of a function and
$\smooth(\R^2, \R^2)$ denotes the set of all infinitely smooth
functions from $\R^2$ to $\R^2$. As in~\cite{LoayzaRomero2025} we note
that, with a slight abuse of notation, the term `outer Riemannian metric' will be
used throughout this paper to refer to a Riemannian metric on $\Diffcn$. Using
such outer Riemannian metrics one expects to immediately obtain shape gradients
which deform the ambient space $\mathbb{R}^2$, and thus inducing a
deformation on the shape. In this paper, we consider the
diffeomorphism group \eqref{Diffcn} as shape space since the shape
space $B_e(S^1,\mathbb{R}^2)$, allowing only smooth shape geometries, is
not sufficient to carry out optimization algorithms for fracture
propagation due to the non-smooth fracture geometry (cf.,
e.g.,~\cite{Suchan2024}).
Moreover, due to their good numerical performance
(cf.~\cite{LoayzaRomero2025,LoayzaWelker2025}), we consider
Sobolev-type Riemannian metrics of order $s=1$, which is a representative
Riemannian metric on $\Diffcn$ (cf., e.g., \cite{Michor2007}).

For the numerical experiments, we use the common single-edge notch
test with either tension or shear loads and compare the results to
established phase field models. A major advantage of performing
fracture simulations based on shape optimization is the lack of
parameters from phase-field approaches, which were originally
introduced by~\cite{Ambrosio1990} and applied to fracture problems in,
e.g.,~\cite{Francfort1998,Bourdin2008,Miehe2010a,Ambati2015,Miehe2010}. The
fracture will be clearly defined by the shape, without the
disadvantage of requiring a sufficiently-fine discretization of the
domain, and without having to interpret the phase field values to
determine the actual fracture path.

The paper is structured as follows: 
In \Cref{sec:background} we describe the fracture propagation problem in variational form~\cite{Francfort1998}, incorporate strain splitting based on \cite{Miehe2010}, illustrate how the fracture problem can be interpreted as an energy minimization problem, and present how this energy minimization problem can be solved using shape optimization algorithms on manifolds. In \Cref{sec:numerical_experiments} we use two common benchmark problems, and numerically investigate multiple different aspects with regards to the shape optimization problem. At the end, we give a summary and conclude our manuscript in \Cref{sec:summary_and_conclusion}.

\section{Fracture propagation model as energy minimization problem} \label{sec:background}

This section introduces the classical brittle, quasi-static fracture propagation problem based on Griffith's criterion~\cite{Griffith1921} and its variational formulation~\cite{Francfort1998}, together with its adaptation to a classical shape optimization setting in \Cref{sec:FracturePropagationObjectiveFunctional}. Furthermore, we explain how strain splitting affects the functional that is to be minimized and how the splitting is achieved in this manuscript in \Cref{sec:ObtainingEigenvectors}. Then, \Cref{sec:shapeOptimizationApproach} introduces the concept of shape optimization on Riemannian manifolds, and establishes the connection between the shape optimization problem in a classical sense (cf., e.g., \cite{Sokolowski1992}) and the optimization problem on a Riemannian manifold. In \Cref{sec:OptimalitySystem}, we calculate the optimality criteria analytically. The last part, \Cref{sec:shapeGradient}, describes the enforcement of the irreversibility constraint, which is an essential part of the fracture propagation problem.

\subsection{Fracture propagation} \label{sec:FracturePropagationObjectiveFunctional}
\begin{figure}[tbp]
	\centering
	\setlength\figureheight{.4\textwidth} 
	\setlength\figurewidth{.4\textwidth}
	\setlength\svgwidth{\figurewidth}
	\begingroup%
  \makeatletter%
  \providecommand\color[2][]{%
    \errmessage{(Inkscape) Color is used for the text in Inkscape, but the package 'color.sty' is not loaded}%
    \renewcommand\color[2][]{}%
  }%
  \providecommand\transparent[1]{%
    \errmessage{(Inkscape) Transparency is used (non-zero) for the text in Inkscape, but the package 'transparent.sty' is not loaded}%
    \renewcommand\transparent[1]{}%
  }%
  \providecommand\rotatebox[2]{#2}%
  \newcommand*\fsize{\dimexpr\f@size pt\relax}%
  \newcommand*\lineheight[1]{\fontsize{\fsize}{#1\fsize}\selectfont}%
  \ifx\svgwidth\undefined%
    \setlength{\unitlength}{317.3126437bp}%
    \ifx\svgscale\undefined%
      \relax%
    \else%
      \setlength{\unitlength}{\unitlength * \real{\svgscale}}%
    \fi%
  \else%
    \setlength{\unitlength}{\svgwidth}%
  \fi%
  \global\let\svgwidth\undefined%
  \global\let\svgscale\undefined%
  \makeatother%
  \begin{picture}(1,0.99725231)%
    \lineheight{1}%
    \setlength\tabcolsep{0pt}%
    \put(0,0){\includegraphics[width=\unitlength,page=1]{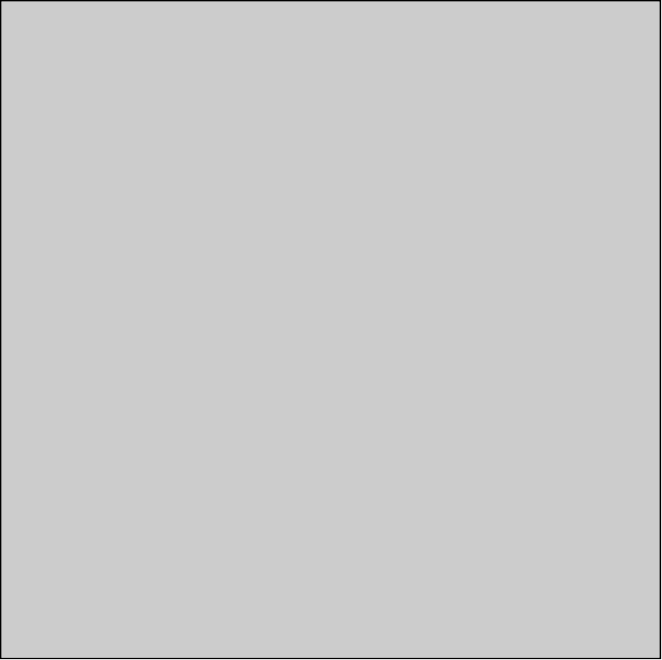}}%
    \put(0.26719752,0.82157014){\color[rgb]{0,0,0}\makebox(0,0)[t]{\lineheight{1.25}\smash{\begin{tabular}[t]{c}$\holdalldomain$\end{tabular}}}}%
    \put(0,0){\includegraphics[width=\unitlength,page=2]{initial_crack_1d_svg-tex.pdf}}%
    \put(0.49989836,0.21126465){\color[rgb]{0,0,0}\makebox(0,0)[t]{\lineheight{1.25}\smash{\begin{tabular}[t]{c}$\physdomain=\holdalldomain \setminus \fracture$\end{tabular}}}}%
    \put(0.72607985,0.52502306){\color[rgb]{0,0,0}\makebox(0,0)[rt]{\lineheight{1.25}\smash{\begin{tabular}[t]{r}$\color{red}\fracture$\end{tabular}}}}%
    \put(0.88083813,0.93269558){\color[rgb]{0,0,0}\makebox(0,0)[rt]{\lineheight{1.25}\smash{\begin{tabular}[t]{r}$\partial\holdalldomain$\end{tabular}}}}%
  \end{picture}%
\endgroup%

	\caption{Sketch of a two-dimensional hold-all domain $\holdalldomain$ that is decomposed into a physical domain~$\physdomain$ and a lower dimensional fracture~$\fracture$, with its boundary $\partial\holdalldomain$.}
	\label{fig:sketchHoldAllDomainAndFracture}
\end{figure}

We consider a hold-all domain $\holdalldomain \subset \R^d$, to be decomposed into a physical domain $\physdomain$
and a lower dimensional fracture $\fracture$, i.e., $\physdomain = \holdalldomain \setminus \fracture$, cf. \Cref{fig:sketchHoldAllDomainAndFracture}. We introduce a time discretization $t_i$, $i=1,\ldots,N$, and an initial fracture $\fracture_0$. Let $\sym\left(\bm{A}\right) = \frac{1}{2}(\bm{A} + \bm{A}^\top)$ be the symmetric part of a matrix, $\strain(\state) = \sym\left(\nabla \state\right) = \frac{1}{2} \left( \nabla \state + \nabla \state^\top \right)$ the linearized strain tensor and $\bm{A} : \bm{B}=\sum_{i=1}^{d} \sum_{j=1}^{d} A_{i j} B_{i j}$ the standard Frobenius product of two matrices. Furthermore, the stress tensor is given as $\stress(\state) = \mathbb{C}:\strain(\state) := 2 \mu \, \strain(\state) + \lambda \trace(\strain(\state)) \, \bm{I}$ with the Lamé parameters $\lambda,\mu > 0$. Let $H^1_D(\physdomain, \R^d)$ be the Hilbert space of vector fields defined on $\physdomain$ whose weak derivatives up to order~$1$ exist, and which additionally have homogenous
Dirichlet boundary data on a prescribed part of the boundary~$\Gamma_D
\subset \partial \physdomain$ with (($d-1$)-dimensional Lebesgue
measure) $|\Gamma_D| \neq 0$. Additionally, forces acting on the volume as well as boundary forces are omitted, $G_c > 0$ denotes the fracture toughness, and $\hmetric^{d-1}$ is the $(d-1)$-dimensional Hausdorff measure. Then, following~\cite{Francfort1998}, a
simple time discrete variational evolution for the displacements $\state^{i+1} \colon \physdomain \rightarrow \R^d$ and fractures $\fracture_{i+1}$ is given by the sequence of $N$ minimization problems
\begin{align}
	\state^{i+1}, \fracture_{i+1} = \argmin_{\state,\fracture} \int_\physdomain \frac{1}{2} \, \stress(\state) : \strain(\state) \d \bm{x} + G_c \hmetric^{d-1}(\fracture)
	\label{eqn:MinimizationWithRespectToStateAndShape}
\end{align}
over all admissible displacements $\state \in H^1_D(\physdomain;\R^d) + \state_D$ and admissible growing fractures $\fracture \supset \fracture_n$, where $\bm{w}_D = \bm{w}_D(t_{i+1})\in H^1_D(\physdomain;\R^d)$ is a given extension of inhomogeneous Dirichlet data into $\physdomain$. With the given $\stress(\state)$, the functional, specifically the first term, can also be expressed as
\begin{align*}
	\int_\physdomain \frac{1}{2} \, \stress(\state) : \strain(\state) \d \bm{x}&= \int_\physdomain \frac{1}{2} \left(2 \mu \strain(\state) + \lambda \trace(\strain(\state)) \bm{I}\right) : \strain(\state) \d \bm{x}\\
	&= \int_\physdomain \mu \, \strain(\state) : \strain(\state) + \frac{\lambda}{2} \trace(\strain(\state))^2 \d \bm{x},
\end{align*}
cf., e.g.,~\cite{Miehe2010}.

As proposed in~\cite{Suchan2023}, assuming sufficient regularity of the domain $\physdomain$,
the minimization with respect to $\state$ can be eliminated to obtain the following minimization problem:
\begin{align}	\label{eqn:fracturebase}
\begin{aligned}
  \min_{\fracture}&\; \int_\physdomain \frac{\lambda}{2} \trace(\strain(\state))^2 + \mu \, \strain(\state) : \strain(\state) \d \bm{x} + G_c \hmetric^{d-1}(\fracture)\\
  &\;\text{s.t. } \int_\physdomain \stress(\state) : \strain(\adjstate) \d \bm{x} = 0  \qquad \forall \adjstate \in H^1_D(\physdomain,\R^d)
\end{aligned}
\end{align}
It should be noted that the first two terms in the objective depend on $\fracture$ by the relation $\physdomain = \holdalldomain \setminus \fracture$. The constraint equation represents the linear elasticity equation in weak form.

According to~\cite{Miehe2010} it is physically unrealistic that compressive strains drive fracture growth in this problem setting. Therefore, it has been proposed to decompose the bulk energy into a tensile and a compressive part. This decomposition is based on the diagonalization of the strain tensor
\begin{align}
\bm{Q}(\state)^{-1} \, \strain(\state) \, \bm{Q}(\state) \eqqcolon \bm{\Sigma}(\state) = \begin{pmatrix}
	\varepsilon_1(\state) & & 0 \\
	 & \ddots & \\
	 0 & & \varepsilon_d(\state)
\end{pmatrix}
\label{eq:StrainDecompositionDiagonalMatrix}
\end{align}
with a matrix $\bm{\Sigma}(\state)$ that contains the eigenvalues (often called principal strains) $\varepsilon_1(\state),\ldots,\varepsilon_d(\state)$ on the diagonal and $\bm{Q}(\state)$ that contains the corresponding eigenvectors (often called principal strain directions). Due to the symmetry of the strain tensor $\strain(\state)$ we are guaranteed that this diagonalization is possible, real eigenvalues are obtained, and the eigenvectors to different eigenvalues are orthogonal. If these eigenvectors are chosen with unit length then we further obtain $\bm{Q}(\state)^{-1} = \bm{Q}(\state)^\top$. The tensile part of the strain tensor is thus chosen as
\begin{align}
\strain^+(\state) = \bm{Q}(\state) \tunderbrace{\begin{pmatrix}
	\max(0,\varepsilon_1(\state)) & & 0 \\
	& \ddots & \\
	0 & & \max(0,\varepsilon_d(\state))
\end{pmatrix}}{\eqqcolon \bm{\Sigma}_{\text{max}}(\state)} \bm{Q}(\state)^\top
\label{eqn:StrainDecompositionOnlyPositiveStrains}
\end{align}
as positive and negative eigenvalues are related to the tensile and compressive part, respectively. Formulated in accordance with~\eqref{eqn:fracturebase},
the full fracture problem then becomes
\begin{align}	\label{eqn:OptProblemnearlyfinal}
	\begin{aligned}
		\min_{\fracture}&\; \int_\physdomain \frac{\lambda}{2} \max(0,\trace(\strain(\state)))^2 + \mu \, \strain^+(\state) : \strain^+(\state) \d \bm{x} + G_c \hmetric^{d-1}(\fracture)\\
		&\;\text{s.t. }\int_\physdomain \stress(\state) : \strain(\adjstate) \d \bm{x} = 0  \qquad \forall \adjstate \in H^1_D(\physdomain,\R^d),%
	\end{aligned}
\end{align}
cf.~\cite{Ambati2015,Miehe2010}.
Using the invariance of the trace operator under a change of basis and noting that $\bm{Q}(\state)^{-1} = \bm{Q}(\state)^\top$ and therefore $\bm{Q}(\state)^\top \bm{Q}(\state) = \bm{I}$, we get
\begin{align}
	\label{eq:trace}
		\strain^+(\state) : \strain^+(\state) = \trace \left( \bm{Q}(\state) \bm{\Sigma}_{\text{max}}(\state)^2 \bm{Q}(\state)^\top \right) = \trace \left( \bm{\Sigma}_{\text{max}}(\state)^2 \right)
\end{align}
yielding
\begin{align}	\label{eqn:OptProblembeforeShape}
	\begin{aligned}
		\min_{\fracture}&\; \int_\physdomain \frac{\lambda}{2} \max(0,\trace(\strain(\state)))^2 + \mu \trace \left( \bm{\Sigma}_{\text{max}}(\state)^2 \right) \d \bm{x} + G_c \hmetric^{d-1}(\fracture) \\
		&\;\text{s.t. }\int_\physdomain \stress(\state) : \strain(\adjstate) \d \bm{x} = 0  \qquad \forall \adjstate \in H^1_D(\physdomain,\R^d).%
	\end{aligned}
\end{align}

\begin{figure}[tbp]
	\centering
	\setlength\figureheight{.4\textwidth} 
	\setlength\figurewidth{.4\textwidth}
	\newlength\svgwidth
	\setlength\svgwidth{\figurewidth}
	\begingroup%
  \makeatletter%
  \providecommand\color[2][]{%
    \errmessage{(Inkscape) Color is used for the text in Inkscape, but the package 'color.sty' is not loaded}%
    \renewcommand\color[2][]{}%
  }%
  \providecommand\transparent[1]{%
    \errmessage{(Inkscape) Transparency is used (non-zero) for the text in Inkscape, but the package 'transparent.sty' is not loaded}%
    \renewcommand\transparent[1]{}%
  }%
  \providecommand\rotatebox[2]{#2}%
  \newcommand*\fsize{\dimexpr\f@size pt\relax}%
  \newcommand*\lineheight[1]{\fontsize{\fsize}{#1\fsize}\selectfont}%
  \ifx\svgwidth\undefined%
    \setlength{\unitlength}{317.31125959bp}%
    \ifx\svgscale\undefined%
      \relax%
    \else%
      \setlength{\unitlength}{\unitlength * \real{\svgscale}}%
    \fi%
  \else%
    \setlength{\unitlength}{\svgwidth}%
  \fi%
  \global\let\svgwidth\undefined%
  \global\let\svgscale\undefined%
  \makeatother%
  \begin{picture}(1,0.99725236)%
    \lineheight{1}%
    \setlength\tabcolsep{0pt}%
    \put(0.26719651,0.82157161){\color[rgb]{0,0,0}\makebox(0,0)[t]{\lineheight{1.25}\smash{\begin{tabular}[t]{c}$\holdalldomain$\end{tabular}}}}%
    \put(0.49989835,0.21126345){\color[rgb]{0,0,0}\makebox(0,0)[t]{\lineheight{1.25}\smash{\begin{tabular}[t]{c}$\physdomain=\holdalldomain \setminus (\remainderdomain \cup \shape)$\end{tabular}}}}%
    \put(0.53945603,0.58684178){\color[rgb]{0,0,0}\makebox(0,0)[t]{\lineheight{1.25}\smash{\begin{tabular}[t]{c}$\color{red}\shape=\partial\remainderdomain$\end{tabular}}}}%
    \put(0.74935388,0.4768989){\color[rgb]{0,0,0}\makebox(0,0)[t]{\lineheight{1.25}\smash{\begin{tabular}[t]{c}$\color{red}\remainderdomain$\end{tabular}}}}%
    \put(0,0){\includegraphics[width=\unitlength,page=1]{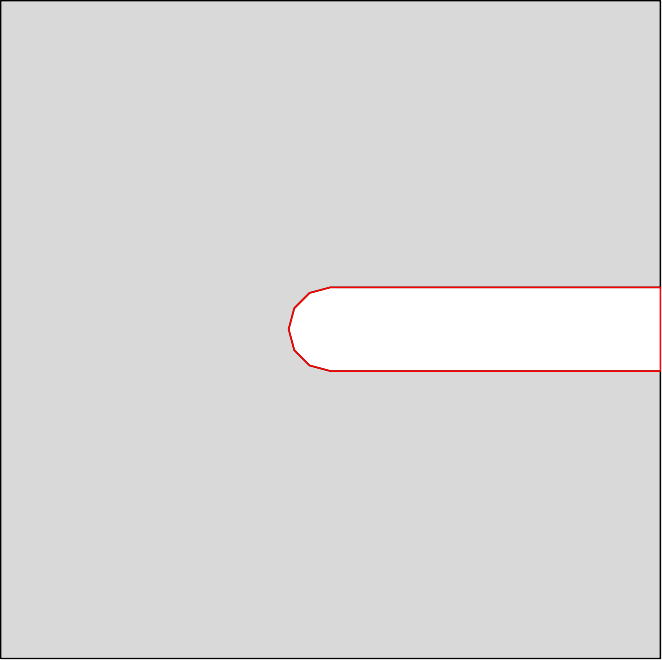}}%
  \end{picture}%
\endgroup%

	\caption{Replacement of the lower dimensional fracture~$\fracture$ (cf. \Cref{fig:sketchHoldAllDomainAndFracture}) by a two-dimensional subdomain~$\remainderdomain$ of the hold-all domain $\holdalldomain$.}
	\label{fig:sketchHoldAllDomainAndRemainderDomain}
\end{figure}

For the proposed shape optimization approach, we replace the lower dimensional fracture
$\fracture$ by a $d$-dimensional subdomain~$\remainderdomain$ of the hold-all domain $\holdalldomain$ and denote the boundary of the domain $\remainderdomain$ as $\shape \coloneqq \partial \remainderdomain$, as sketched in \Cref{fig:sketchHoldAllDomainAndRemainderDomain}. Thus, we have $\physdomain = \holdalldomain \setminus \left(\remainderdomain \cup \shape\right)$.  The subset~$\remainderdomain$ now describes a volume no longer attached to
the physical domain $\physdomain$, i.e., $\holdalldomain \supset \shape
\cup \physdomain$
 and a spurious remainder $\remainderdomain$, of which we need to disincentivize the increase in volume, as a brittle fracture should not remove any bulk material. We achieve this by adding an additional term to the objective functional with $\nu > 0$ that increases when the volume of $\remainderdomain$ is increasing, e.g., $\int_{\remainderdomain} \nu \d \bm{x}$.
Consequently, we model irreversibility of the fracture by requiring that $\left( \remainderdomain_{i} \cup \shape_{i} \right) \subset \left( \remainderdomain_{i+1} \cup \shape_{i+1} \right)$ at each time step $t_i$, $i=1,\ldots,N$, or equivalently
$
\physdomain_{i+1} \subset \physdomain_{i}.
$

As the width of the fracture can and will be chosen to be negligible in comparison to the length, the length of $\shape$ is now, approximately, twice the length of the fracture. Furthermore, adding a constant value (here: two times the width of the fracture) to the length of the fracture and therefore to the objective functional does not influence the solution of the optimization problem. Thus, the $(d-1)$-dimensional Hausdorff measure of $\mathcal{F}$ is approximated as half of the length of $\shape$ in order to rescale the objective accordingly. We further observe that the additional term in the objective functional (that disincentivizes an increase of volume of $\remainderdomain$) can be expressed as
\begin{align*}
	\int_{\remainderdomain} \nu \d \bm{x} = \int_{\holdalldomain\setminus(\physdomain \cup \shape)} \nu \d \bm{x} = \int_{\holdalldomain} \nu \d \bm{x} - \int_{\physdomain} \nu \d \bm{x} - \int_{\shape} \nu \d \bm{x}
\end{align*}
and, since the integral over $\holdalldomain$ never changes (and therefore also only constitutes a constant added to the objective functional) and the integral over $\shape$ is zero, the final optimization problem for $i=1,\ldots,N$ considered in this manuscript reads
\begin{subequations}
\label{eqn:OptProblemFinal}
\begin{align}	\label{eqn:OptProblemFinal_1}
		\min_{\remainderdomain}&\; \tunderbrace{\tunderbrace{\int_\physdomain \frac{\lambda}{2} \max(0,\trace(\strain(\state)))^2 + \mu \trace \left( \bm{\Sigma}_{\text{max}}(\state)^2 \right) \d \bm{x}}{\eqqcolon E_{\text{bulk}}} + \tunderbrace{G_c \int_{\shape} \frac{1}{2} \d \bm{s}}{\eqqcolon E_{\text{frac}}} - \tunderbrace{\int_{\physdomain} \nu \d \bm{x}}{\eqqcolon E_{\text{reg}}}}{\eqqcolon J(\remainderdomain)}\\
		\label{eqn:OptProblemFinal_2}
		&\;\text{s.t. }\tunderbrace{\int_\physdomain \stress(\state) : \strain(\adjstate) \d \bm{x}}{\eqqcolon c(\remainderdomain,\state,\adjstate)} = 0 \qquad \forall \adjstate \in H^1_D(\physdomain,\R^d), \\
		\label{eqn:OptProblemFinal_3}
		&\;\text{\phantom{s.t. }}\physdomain \subset \physdomain_{i}.
\end{align}
\end{subequations}

\subsection{Obtaining the eigenvectors} \label{sec:ObtainingEigenvectors}

One possible way to obtain the matrix $\bm{Q}(\state)$ is to perform a standard eigenvalue-eigenvector decomposition. For this manuscript, we consider $d=2$, for which a general symmetric strain tensor is given by
\begin{align}
	\strain(\state) = \begin{pmatrix}
		\varepsilon_{1 1}(\state) & \varepsilon_{1 2}(\state) \\
		\varepsilon_{2 1}(\state) & \varepsilon_{2 2}(\state)
	\end{pmatrix} = \frac{1}{2} \left( \nabla \state + {\nabla \state}^\top \right).
	\label{eq:strainTensorDefinition}
\end{align}
As $\strain(\state)$ is a symmetric tensor, we have $\varepsilon_{1 2}(\state)=\varepsilon_{2 1}(\state)$.
If $\varepsilon_{1 2}(\state) = 0$, then $\strain(\state)$ is a diagonal matrix. Thus, if $\varepsilon_{1 2}(\state) = 0$, then the eigenvalues of $\strain(\state)$, which we denote as $\varepsilon_{1}(\state)$ and $\varepsilon_{2}(\state)$, are its diagonal entries and the (normalized) eigenvectors are $\bm{q}_{1}(\state) = \left( 1, 0 \right)^\top$ and $\bm{q}_{2}(\state) = \left( 0, 1 \right)^\top$.
If $\varepsilon_{1 2}(\state) \neq 0$, then the eigenvalues of $\strain(\state)$ can be determined by solving $\det(\strain(\state) - \varepsilon_i(\state) \bm{I}) \hastobe 0$ for $\varepsilon_i(\state)$, $i=1,2$, and are given by
\begin{align*}
	\varepsilon_{1,2}(\state) &= \frac{\trace(\strain(\state))}{2} \pm \sqrt{\frac{\trace(\strain(\state))^2}{4} - \det(\strain(\state))} \\
	&= \frac{\trace(\strain(\state))}{2} \pm \sqrt{-\frac{\trace(\strain(\state))^2}{4} + \frac{\trace\left(\strain(\state)^2\right)}{2}}.
\end{align*}
The eigenvectors for $\varepsilon_{1 2}(\state) \neq 0$ can be calculated, for which we omit the calculation here, and result in
\begin{align*}
  \bm{q}_{1,2}(\state) = \begin{pmatrix}
    1 \\
    \frac{\varepsilon_{2 2}(\state) - \varepsilon_{1 1}(\state)}{2 \varepsilon_{1 2}(\state)} \pm \sqrt{\frac{\left(\varepsilon_{2 2}(\state) - \varepsilon_{1 1}(\state)\right)^2}{4 \varepsilon^2_{1 2}(\state)} + 1}
  \end{pmatrix}.
\end{align*}
The additional normalization, which is required such that $\bm{Q}(\state)^{-1} = \bm{Q}(\state)^\top$ holds (cf.~\eqref{eqn:StrainDecompositionOnlyPositiveStrains}), would further complicate this expression.
An alternative, possibly simpler approach is to obtain the eigenvectors by observing that $\bm{Q}(\state)$ is a rotation matrix, as commonly encountered during the planar stress-state investigations using Mohr's circle, cf., e.g.,~\cite{Mohr1906,Gross2014}. For $d=2$ the rotation is described by one rotation angle $\alpha(\state)$, and the rotation matrix then is given as
\begin{align}
	\label{eq:RotationMatrixQ}
	\bm{Q}(\state) = \begin{pmatrix}
		\cos(\alpha(\state)) & -\sin(\alpha(\state)) \\
		\sin(\alpha(\state)) & \cos(\alpha(\state))
	\end{pmatrix}.
\end{align}
By setting the expressions on the off-diagonal of
$\bm{Q}(\state)^\top \, \strain(\state) \, \bm{Q}(\state)$ to $0$ and
solving for the angle, it can be obtained that $\alpha(\state) =
\frac{1}{2} \arctan\left(\frac{2 \varepsilon_{1
      2}(\state)}{\varepsilon_{1 1}(\state)-\varepsilon_{2
      2}(\state)}\right)$. For later derivations, let $\xi(\state) \coloneqq \frac{2 \varepsilon_{1
      2}(\state)}{\varepsilon_{1 1}(\state)-\varepsilon_{2
      2}(\state)}$ such that $\alpha(\state) = \frac{1}{2} \arctan\left(\xi(\state)\right)$.
We observe that only for $\varepsilon_{1 1}(\state) = \varepsilon_{2 2}(\state)$ and at the same time  $\varepsilon_{1 2}(\state)=0$ this expression for $\alpha(\state)$ is possibly undefined. However, we also observe that this case already yields a diagonal $\strain(\state)$, since in this case we have
\begin{align*}
	\strain(\state) = \begin{pmatrix}
		\varepsilon_{1 1}(\state) & 0 \\
		0 & \varepsilon_{1 1}(\state)
	\end{pmatrix},
\end{align*}
and therefore the rotation matrix is given by the identity matrix, i.e., $\bm{Q}(\state) = \bm{I}$.
In order to have a consistent order of eigenvalues
$\varepsilon_{1}(\state) \leq \varepsilon_{2}(\state)$, and by
considering the definition of the function $\operatorname{arctan2}$ as described in, e.g., \cite[Section 7.4.1]{Hughes2014},
we express the rotation angle $\alpha(\state)$ as
\begin{align}
  \label{eq:RotationAngleAlpha}
  \alpha(\state) = \frac{1}{2} \arctan\left(\frac{2 \varepsilon_{1 2}(\state)}{\varepsilon_{1 1}(\state)-\varepsilon_{2 2}(\state)}\right) + \begin{cases}
    0 & \text{, if } \varepsilon_{1 1}(\state) \leq \varepsilon_{2 2}(\state) \\
    \frac{\pi}{2} & \text{, else.}
  \end{cases}
\end{align}
We provide a numerical demonstration of this eigenvalue decomposition to obtain the strain splitting in~\Cref{app:demonstration_strain_splitting} for the interested reader, even though this is not the main focus of the manuscript.

\subsection{Shape optimization approach based on manifolds} \label{sec:shapeOptimizationApproach}

In order to facilitate the modeling of the optimization problem~\eqref{eqn:OptProblemFinal} as a problem on a Riemannian manifold, we restrict ourselves to $d=2$ for the rest of this manuscript. 

Recently, fracture propagation has already been considered in the context of shape optimization on manifolds (cf.~\cite{Suchan2024,Suchan2023}).  In~\cite{Suchan2023}, the fracture is modeled by a smooth curve, i.e., kinks in the fracture geometry are not allowed.
In order to avoid this limitation,~\cite{Suchan2024} models a fracture by a piecewise-smooth curve and defines a corresponding space involving these fractures based on the product manifold established in~\cite{Pryymak2023}. 
Although both publications have different shape spaces and, therefore, perform optimization on different spaces, they consider the same Riemannian metric on the corresponding spaces involving the fractures, which is the Steklov-Poincaré Riemannian metric defined in~\cite{Schulz2016a}.
In this paper, we will concentrate on so-called outer Riemannian metrics, which is the main contrast to recent publications~\cite{Suchan2024,Suchan2023}. 
One reason why we consider outer Riemannian metrics in this paper is the overall better performance in comparison to inner Riemannian metrics like the Steklov-Poincaré Riemannian metric,
cf.,~\cite{LoayzaRomero2025,LoayzaWelker2025}. 
In particular, the $s$-th order Sobolev type Riemannian metrics, which we consider in this paper, have nice theoretical properties under special conditions like geodesic completeness and hence the global existence of the exponential map (cf.~\cite{Bauer2013overview}), which are helpful for, e.g., convergence proofs, computation of next iterates during an optimization procedure, etc.
The use of outer Riemannian metrics---in particular Sobolev-type Riemannian metrics on the diffeomorphism group---in the context of PDE-constrained shape optimization was first studied in~\cite{LoayzaRomero2025}.
The main advantage of using outer Riemannian metrics is that we obtain shape gradients which deform the ambient space. In contrast, inner Riemannian metrics would only deform the fracture itself and 
not the ambient space, and as we need to solve PDE, e.g., the state equation~\eqref{eqn:OptProblemFinal_2}, on the ambient space, thus the required discretization of the ambient space quickly deteriorates if only the fracture were deformed, instead of the fracture together with the ambient space. The topic of extension operators for shape gradients obtained using inner Riemannian metrics to the ambient space has been considered in several publications, cf., e.g., \cite{Gangl2015,Schulz2016,Onyshkevych2021}, but can be avoided altogether using outer Riemannian metrics.

\begin{figure}[tbp]
	\centering
	\setlength\figureheight{.4\textwidth} 
	\setlength\figurewidth{.95\textwidth}
	\newlength\svgwidth
	\setlength\svgwidth{\figurewidth}
	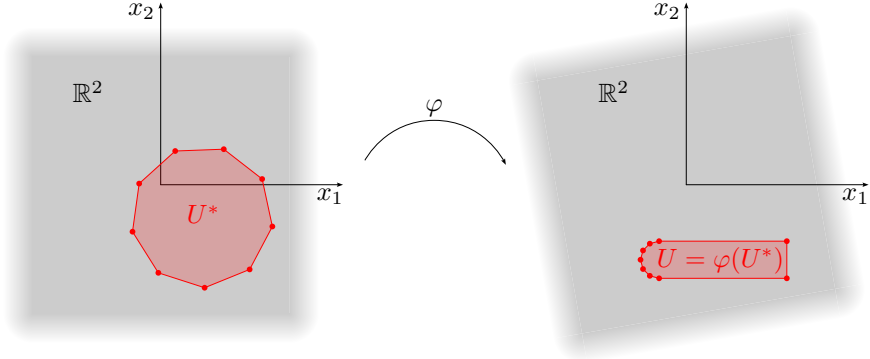
	\caption{Sketch of one exemplary diffeomorphism $\varphi$ from $\R^2$ to $\R^2$. In this example, the interior of a convex regular polygon~$\remainderdomain^*$, a subset of $\R^2$, is mapped to~$\remainderdomain$.}
	\label{fig:DiffeomorphismR2ToR2}
\end{figure}

In this paper, inspired by~\cite{LoayzaRomero2025}, we consider the set~\eqref{Diffcn} as the shape space for~\eqref{eqn:OptProblemFinal}. 
Therefore, the domain $\remainderdomain$ from \Cref{sec:FracturePropagationObjectiveFunctional} can be described as the image of a diffeomorphism $\varphi \in \Diffcn$ applied to a non-empty measurable set $\remainderdomain^* \subset \R^2$, i.e., $\remainderdomain = \varphi(\remainderdomain^*)$. Thus, we will identify $\varphi\in \Diffcn$ with the image of $\remainderdomain^*$. In order to allow kinks in $\remainderdomain$, we consider $\remainderdomain^*$ as the interior of a unit convex regular polygon, cf. \Cref{fig:DiffeomorphismR2ToR2}.
In order to choose the optimal amount $k\in\mathbb{N}$ of kinks in $\remainderdomain$, we assume that an $\epsilon>0$ is chosen such that $k$ is the minimal number such that $\left| A^{\text{poly}}-A^{\text{circle}} \right| \leq \epsilon$, where $A^{\text{poly}}$ and $A^{\text{circle}}$ are the enclosed areas of the unit convex regular polygon with $k$ kinks and the unit circle, respectively.
The diffeomorphism group~\eqref{Diffcn}  is a smooth submanifold of $C^\infty(\R^2, \R^2)$ and a regular Lie group with Lie algebra\footnote{For an introduction into Lie theory including Lie groups and Lie algebras we refer to the literature, e.g., \cite{Lee2009,Lie}.} given by the space of diffeomorphisms with compact support (cf.~\cite[Theorem 43.1]{Kriegl1997}).
In order to establish a distance measure on the shape space for gradient-based optimization of an objective functional $j \colon \Diffcn \to \R$, we equip it with a Sobolev-type Riemannian metric of order $s$, which for two elements $\vV,\vW\colon \R^2\to\R^2 $, $\vV=\vV(x_1,x_2)$, $\vW=\vW(x_1,x_2)$, of its tangent space at a point $\varphi \in \Diffcn$, denoted by $T_\varphi \Diffcn$, is defined as
\begin{equation}
	\label{eq:metric1}
	g_\varphi
	^s(\vV,\vW)= \int_{\R^2} \left< (\operatorname{id}-A\Delta)^s\vV,\vW \right> \d \bm{x}\qquad (A>0),
\end{equation}
where $\Delta$ denotes the Laplace operator (cf.~\cite{Bauer2013overview,Bauer2020sobolev}). 
As usual, we define a shape functional on $(\Diffcn, g^s)$ as a function $j\colon \Diffcn \to \R$ and denote the space of all shape functionals on $\Diffcn$ by $\mathcal{A}(\Diffcn)$. 
Then, the Riemannian gradient of $j$ at $\varphi$ is the element~$\vV \in T_\varphi \Diffcn$ that fulfills
\begin{align}
	\label{eqn:ShapeGradientPushforward}
	g_\varphi^s(\vV,\vW) = \mathfrak{d} j(\varphi) [\vW] \quad\forall \vW\in T_\varphi \Diffcn,
\end{align}
where the right-hand side denotes the pushforward of $j$ at $\varphi$ in direction $\vW$. However, we refrain from providing a thorough introduction and definition of the pushforward in this manuscript due to its more applied focus and a theorem in~\cite{LoayzaRomero2025}.\footnote{ The theorem in~\cite{LoayzaRomero2025} provides the connection of the pushforward as a differential of a $j$, where $j$ is is seen as a function between two manifolds, and the differential in classical shape optimization. The interested reader is referred to the literature, e.g.,~\cite{Lee2009}, for a complete mathematical definition of the pushforward. %
If we consider
\begin{align}
	\label{eqn:ObjectiveFunctionalOnShapeSpace}
	\mathcal{J} \coloneqq\{j \in \mathcal{A}(\Diffcn) \colon \exists J \colon \mathcal{D}  \subset P(\R^2) \to \R 
	\text{ such that }j(\varphi)=J(\varphi(\remainderdomain^*)) \ \forall \varphi \in \Diffcn \},
\end{align} 
where $P(\R^2)$ denotes the power set of $\R^2$,
then any  $j \in \mathcal{J}$ can be represented by a shape functional $J$ from classical shape optimization theory (cf.~\cite{LoayzaRomero2025}). This permits the mentioned connection of the pushforward and the differential in classical shape optimization, which we use in \Cref{sec:OptimalitySystem}.}%

It is worth mentioning that since we are only interested in the fracture itself, i.e., the image of the diffeomorphism applied to $\remainderdomain^*$, in general we need to consider the quotient space
$\Diffcn/\Diffc(\R^2, \remainderdomain^*)$, where $ \Diffc(\R^2, \remainderdomain^*)\coloneqq \{ \varphi \in \Diffcn \colon \varphi( \remainderdomain^*)=\remainderdomain^*\}$. 
Please note that this quotient space is only necessary to avoid that the images of the diffeomorphisms (representing fractures), which only differ by reparametrizations, are considered as equal to each other. For simplicity reasons, we only consider~\eqref{Diffcn} as the shape space and not the quotient space since the value of the objective functional in the optimization problem will not change under reparametrizations.

Now, we will describe how the optimization problem~\eqref{eqn:OptProblemFinal} fits into the setting that was just described.
In order to achieve this, as in~\cite{LoayzaRomero2025}
we now want to only consider elements in $\Diffcn$ which are the identity outside of some bounded domain. For our example, the bounded domain in the image of $\Diffcn$ corresponds to $\holdalldomain$. We will denote the elements in $\Diffcn$ that fulfill this condition by $\Diffc^\holdalldomain(\R^2)$. Furthermore, let $\holdalldomain^* \subset \R^2$ be a reference domain in the preimage of $\Diffcn$, and define $\remainderdomain^* \subset \holdalldomain^*$ and $\physdomain^* = \holdalldomain^* \setminus \left( \remainderdomain^* \cup \partial \remainderdomain^* \right)$ such that $\holdalldomain = \varphi(\holdalldomain^*)$, $\remainderdomain = \varphi(\remainderdomain^*)$ and $\physdomain = \varphi(\physdomain^*)$. Due to~\eqref{eqn:ObjectiveFunctionalOnShapeSpace}, the optimization problem analogous to~\eqref{eqn:OptProblemFinal} using $\varphi$ is then
\begin{subequations}
	\label{eqn:OptProblemFinalShapeSpace}
	\begin{align}	\label{eqn:OptProblemFinalShapeSpace_1}
			&\begin{aligned}
			\min_{\varphi \in \Diffc^\holdalldomain(\R^2)}&\; \bigg( \tunderbrace{\int_{\varphi(\physdomain^*)} \frac{\lambda}{2} \max(0,\trace(\strain(\state)))^2 + \mu \trace \left( \bm{\Sigma}_{\text{max}}(\state)^2 \right) \d \bm{x}}{\eqqcolon E_{\text{bulk}}}  \\
			&\hphantom{\bigg(\;}\tunderbrace{+ \tunderbrace{G_c \int_{\varphi(\partial\remainderdomain^*)} \frac{1}{2} \d \bm{s}}{\eqqcolon E_{\text{frac}}} - \tunderbrace{\int_{\varphi(\physdomain^*)} \nu \d \bm{x}}{\eqqcolon E_{\text{reg}}} \bigg) \hspace*{2.11cm} }{j(\varphi)}
			\end{aligned}
			\\
			\label{eqn:OptProblemFinalShapeSpace_2}
			&\hspace*{2em}\;\text{s.t. }\tunderbrace{\int_{\varphi(\physdomain^*)} \stress(\state) : \strain(\adjstate) \d \bm{x}}{\eqqcolon c(\varphi,\state,\adjstate)} = 0 \qquad \forall \adjstate \in H^1_D({\varphi(\physdomain^*)},\R^2), \\
			\label{eqn:OptProblemFinalShapeSpace_3}
			&\hspace*{2em}\;\text{\phantom{s.t. }}{\varphi(\physdomain^*)} \subset \physdomain_{i}.
	\end{align}
\end{subequations}
Correspondingly, tangent vectors are elements in $C_0^\infty(\holdalldomain,\R^2)$ (for details we refer to~\cite{LoayzaRomero2025}), where $C_0^\infty(\holdalldomain,\R^2)$ denotes the set of all infinitely smooth functions mapping from $\holdalldomain$ to $\R^2$ which vanish on $\partial \holdalldomain$. %
In order to optimize on the shape space $\Diffc^\holdalldomain(\R^2)$ based on gradient information we not only require a Riemannian metric but also a method of mapping from the tangent space of the manifold to the manifold itself in order to formulate the next (shape-)iterate in an algorithm. The  exponential map is the first option, however due to the numerical effort of solving a second-order differential equation for the exponential map, it is common to use a retraction instead. In particular, in~\cite{LoayzaWelker2025}, it is shown in an experiment that even without the use of geodesics, outer Riemannian metrics have overall better performance than, e.g., the inner Steklov-Poincaré Riemannian metric. We focus on the same retraction on the diffeomorphism group as in~\cite[Subsection 5.1.1]{LoayzaRomero2025}. This retraction is only a first-order local approximation of the exponential map, and thus it will be well defined only for small enough updates in the direction of the tangent vectors. Since this retraction can be interpreted as a version of the perturbation of identity method, which is well-known in the context of shape optimization, %
we can define a domain update by
\begin{align}
	\label{eqn:DomainUpdate}
	\holdalldomain_\tau = \{ \bm{x} \in \holdalldomain \colon \bm{x} + \tau \vV (\bm{x}) \}
\end{align}
with a vector field $\vV \in C_0^\infty(\holdalldomain,\R^2)$ and $\tau \geq 0$, which in turn also updates $\remainderdomain$.

\subsection{Optimality system} \label{sec:OptimalitySystem}

As already described in \Cref{sec:shapeOptimizationApproach}, we have a correspondence of~\eqref{eqn:OptProblemFinal} and~\eqref{eqn:OptProblemFinalShapeSpace}. Moreover, we consider the set~\eqref{eqn:ObjectiveFunctionalOnShapeSpace} which connects the pushforward with the derivative in classical shape optimization thanks to~\cite[Theorem~3.1]{LoayzaRomero2025}. This enables us to consider~\eqref{eqn:OptProblemFinal} and use classical (shape) optimization approaches.
In order to define the optimality condition of the system~\eqref{eqn:OptProblemFinal_1}--\eqref{eqn:OptProblemFinal_2}, where the minimization is performed with respect to $\remainderdomain$, we first need to define the Lagrange functional, which is given as
\begin{align}
	\label{eq:LagrangeFunction}
	\begin{aligned}
	L(\remainderdomain) &= L(\remainderdomain, \state(\remainderdomain), \adjstate(\remainderdomain)) \\
	&= E_{\text{bulk}}(\remainderdomain,\state(\remainderdomain)) + E_{\text{frac}}(\remainderdomain) + E_{\text{reg}}(\remainderdomain) + c(\remainderdomain,\state(\remainderdomain),\adjstate(\remainderdomain)).
	\end{aligned}
\end{align}

To formulate the necessary optimality conditions, the total (Lagrangian) derivative with respect to $\remainderdomain$ (cf.~\cite[p. 38--39]{Paterson1997} for an explanation of the Lagrangian derivative), denoted by~$\deriv_\remainderdomain$ in the following, of $L$ in the direction of a sufficiently smooth vector field $\vV$ has to be computed. In the following, we denote the partial (Eulerian) derivatives with respect to $\remainderdomain$, $\state$ and $\adjstate$ as $\partial_{\remainderdomain}$, $\partial_{\state}$ and $\partial_{\adjstate}$, respectively. Provided sufficient regularity, the chain rule directly yields
\begin{align}
	\label{eqn:KKTconditions}
	\begin{aligned}
		&\deriv_{\remainderdomain} L(\remainderdomain, \state(\remainderdomain), \adjstate(\remainderdomain)) [\vV] \\
		&= \partial_{\remainderdomain} \left( E_{\text{bulk}}(\remainderdomain, \state) + E_{\text{frac}}(\remainderdomain) + E_{\text{reg}}(\remainderdomain) + c(\remainderdomain, \state, \adjstate) \right) [\vV] \\
		&\hphantom{=\ }+ \partial_{\state} \left( E_{\text{bulk}}(\remainderdomain,\state) + c(\remainderdomain, \state, \adjstate) \right) \, \deriv_{\remainderdomain} \state(\remainderdomain) [\vV] \\
		&\hphantom{=\ }+ \partial_{\adjstate} \left( c(\remainderdomain, \state, \adjstate) \right) \, \deriv_{\remainderdomain} \adjstate(\remainderdomain)[\vV],
	\end{aligned}
\end{align}
where we have already omitted terms with zero derivative. When at
an optimum, these partial derivatives all have to be zero.
Therefore, from setting these partial derivatives to zero we obtain the design equation, the adjoint equation and the state equation. With this, a point $(\hat{\remainderdomain}, \hat{\state}, \hat{\adjstate}) = (\hat{\remainderdomain}, \state(\hat{\remainderdomain}), \adjstate(\hat{\remainderdomain}))$ that fulfills the first-order necessary conditions (sometimes also called a Karush-Kuhn-Tucker (KKT) point~\cite{Karush1939,Kuhn1951}) of minimizing $L(\remainderdomain)$ has to satisfy that~\eqref{eqn:KKTconditions} is equal to zero for all $\vV$.

\paragraph{State equation}
The state equation follows from requiring $0=\partial_{\adjstate} \left( c(\remainderdomain,\state,\adjstate) \right) [\tilde{\adjstate}]$ for all $\tilde{\adjstate}\in H^1_D(\physdomain;\R^2)$ and reads: Find a $\state \in H^1_D(\physdomain;\R^2)+ \state_D$, where $\state_D$ denotes the extension of inhomogeneous Dirichlet boundary data as described in \Cref{sec:FracturePropagationObjectiveFunctional}, such that
\begin{align}
	0 = \int_\physdomain \stress(\state) : \strain(\tilde{\adjstate}) \d \bm{x} \qquad \forall \tilde{\adjstate} \in H^1_D(\physdomain;\R^2).
	\label{eq:StateEquation}
\end{align}
As described in~\Cref{sec:FracturePropagationObjectiveFunctional} we have omitted volume and surface loads. Thus, they do not appear in the state equation. The short calculation can be found in \Cref{app:calculation_state_equation}.

\paragraph{Adjoint equation}
The adjoint equation follows from requiring
\begin{align*}
	0=\partial_{\state} \left( E_{\text{bulk}}(\remainderdomain,\state) + c(\remainderdomain, \state, \adjstate) \right) [\tilde{\state}]
\end{align*}
for all $\tilde{\state}\in H^1_D(\physdomain;\R^2)$. It involves the derivative of the bulk energy~$E_{\text{bulk}}$ with respect to the state~$\state$ and is therefore slightly more challenging than the calculation of the state equation.  It is useful to define the condition
\begin{align}
	\label{eqn:ConditionStrainDecomposition}
	\varepsilon_{1 1}(\state) = \varepsilon_{2 2}(\state) \land \varepsilon_{1 2}(\state)=0
\end{align}
for the following.
The whole calculation can be found in \Cref{app:calculation_adjoint_equation}, which yields: Find $\adjstate \in H^1_D(\physdomain;\R^2)$ such that
\begin{align}
	\label{eq:AdjointEquation}
	\begin{aligned}
		0 &= \int_\physdomain \strain(\tilde{\state}) : \stress(\adjstate) + \lambda \max(0,\trace(\strain(\state))) \cdot \trace(\strain(\tilde{\state})) \\
		&\hphantom{= \int_\physdomain\,}+ 2 \mu \begin{cases}
			\trace \left(\bm{\Sigma}_{\text{max}}(\state) \, \strain(\tilde{\state})\right) \d \bm{x} & \!\!\!\!\text{, if \eqref{eqn:ConditionStrainDecomposition}} , \\
				\frac{\varepsilon_{1 2}(\tilde{\state}) \cdot \left( \varepsilon_{1 1}(\state)-\varepsilon_{2 2}(\state) \right) - \varepsilon_{1 2}(\state) \cdot \left(\varepsilon_{1 1}(\tilde{\state})-\varepsilon_{2 2}(\tilde{\state})\right)}{\left(\varepsilon_{1 1}(\state)-\varepsilon_{2 2}(\state)\right)^2 + 4 \varepsilon^2_{1 2}(\state)} \\
				\ \ \cdot \trace \left( \bm{\Sigma}_{\text{max}}(\state) \left( \bm{R}(\state)^\top \strain(\state) \, \bm{Q}(\state) + \bm{Q}(\state)^\top \strain(\state) \, \bm{R}(\state) \right) \right) \\
				\ + \trace \left( \bm{\Sigma}_{\text{max}}(\state) \, \bm{Q}(\state)^\top \strain(\tilde{\state}) \, \bm{Q}(\state) \right) \d \bm{x}
			 & \!\!\!\!\text{, else}
		\end{cases}
	\end{aligned}
\end{align}
holds for all $\tilde{\state} \in H^1_D(\physdomain;\R^2)$.

\paragraph{Design equation}
As we aim to perform optimization, besides the state and adjoint equation we lastly require the design equation
\begin{align*}
	0=\partial_{\remainderdomain} \left( E_{\text{bulk}}(\remainderdomain, \state) + E_{\text{frac}}(\remainderdomain) + E_{\text{reg}}(\remainderdomain) + c(\remainderdomain, \state, \adjstate) \right)[\vV]
\end{align*}
to be fulfilled for all $\vV \in C_0^\infty(\holdalldomain,\R^2)$. Let a perturbation of $\remainderdomain$ in direction $\vV$ (where $\vV$ is restricted to $\remainderdomain \subset \holdalldomain$) be defined as in~\eqref{eqn:DomainUpdate}.
In general, for an arbitrary differentiable functional $J$, the expression $\shapeDeriv J(\remainderdomain)[\vV]$ is often also called the shape derivative of $J$ in direction of a sufficiently smooth vector field $\vV$, which is defined as the Eulerian derivative
\begin{align}
	\label{eqn:shapeDerivative}
	\shapeDeriv J(\remainderdomain)[\vV] \coloneqq \lim_{\tau \to 0^+} \frac{J(\remainderdomain_\tau) - J(\remainderdomain)}{\tau}
\end{align}
if this derivative exists for all $\vV$ and the mapping $\vV \mapsto \shapeDeriv J(\remainderdomain)[\vV]$ is linear and continuous. %
According to~\cite[Theorem~3.1]{LoayzaRomero2025} and with~\eqref{eqn:ObjectiveFunctionalOnShapeSpace}, the shape gradient (that can later be used to perform optimization) can be obtained from the shape derivative with respect to~$\remainderdomain$ in direction $\vV$, instead of the pushforward, which is why we omitted the definition of the pushforward as it is not strictly needed.
The shape derivative of $L(\remainderdomain)$ in~\eqref{eq:LagrangeFunction} reads
\begin{align}
  \label{eqn:ShapeDerivative}
  \begin{aligned}
    &\shapeDeriv L(\remainderdomain, \state, \adjstate)[\vV] \\
    &= - \int_\physdomain \lambda \max(0,\trace(\strain(\state))) \cdot \trace \left( \nabla \state \nabla \vV \right) \\
    &\hphantom{=- \int_\physdomain\,} + \sym\left( \nabla \state \nabla \vV \right) : \stress(\adjstate) + \stress(\state) : \sym\left( \nabla \adjstate \nabla \vV \right) \d \bm{x} \\
    &\hphantom{=\ }- \int_\physdomain 2 \mu \begin{cases}
      \trace \left( \bm{\Sigma}_{\text{max}}(\state) \cdot  \sym\left( \nabla \state \nabla \vV \right) \right) \d \bm{x} & \!\!\!\!\text{, if \eqref{eqn:ConditionStrainDecomposition}}, \\[10pt]
        \frac{ \left( \begin{pmatrix}
              -\varepsilon_{1 2}(\state) & \frac{\varepsilon_{1 1}(\state)-\varepsilon_{2 2}(\state)}{2} \\
              \frac{\varepsilon_{1 1}(\state)-\varepsilon_{2 2}(\state)}{2} & \varepsilon_{1 2}(\state)
            \end{pmatrix} {\nabla \state}^\top \right) : \nabla \vV}{\left(\varepsilon_{1 1}(\state)-\varepsilon_{2 2}(\state)\right)^2 + 4 \varepsilon^2_{1 2}(\state)} \\
        \ \ \cdot \trace \left( \bm{\Sigma}_{\text{max}}(\state) \cdot \left( \bm{R}(\state)^\top \strain(\state) \, \bm{Q}(\state) + \bm{Q}(\state)^\top \strain(\state) \, \bm{R}(\state) \right) \right) \\
        \ + \trace \left( \bm{\Sigma}_{\text{max}}(\state) \cdot \bm{Q}(\state)^\top \sym\left( \nabla \state \nabla \vV \right) \, \bm{Q}(\state) \right) \d \bm{x}  & \!\!\!\!\text{, else}
    \end{cases} \\
    &\hphantom{=\ }+ \int_\physdomain \Div(\vV) \cdot \left( \frac{\lambda}{2} \max(0,\trace(\strain(\state)))^2 + \mu \trace \left( \bm{\Sigma}_{\text{max}}(\state)^2 \right) + \stress(\state) : \strain(\adjstate) - \nu \right) \d \bm{x} \\
    &\hphantom{=\ }+ \int_\shape \frac{G_c}{2} \left(\Div(\vV) - \bm{n}^\top \nabla \vV \bm{n}\right) \d \bm{s}.
  \end{aligned}\hspace*{-1em}
\end{align}
The detailed calculation can be found in \Cref{app:calculation_shape_derivative}. Setting the shape derivative to zero for all $\vV$ gives a necessary condition for a minimum of the optimization problem~\eqref{eqn:OptProblemFinal}.

\subsection{Shape gradient with enforced irreversibility constraint} \label{sec:shapeGradient}
As already mentioned in \Cref{sec:OptimalitySystem}, the shape gradient is obtained from the shape derivative with respect to~$\remainderdomain$ in direction $\vV$, which was calculated in the previous section, instead of the pushforward.
We have in mind that
we want to only consider elements in 
$\Diffcn$ which are the identity outside of some bounded domain $\holdalldomain\subset\R^2$, i.e., $\Diffc^\holdalldomain(\R^2)$ as already mentioned in \Cref{sec:shapeOptimizationApproach}. Here, the tangent
vectors are elements in $C_0^\infty(\holdalldomain,\R^2)$.
Then, for $L(\remainderdomain)$ to be minimized as in~\eqref{eq:LagrangeFunction}, the
Riemannian gradients $\vV = \nabla_R L(\remainderdomain)$
associated with the Riemannian metrics $g^s$ at $\varphi$ defined in~\eqref{eq:metric1} can be computed by solving
\begin{align}
	\label{eq:shape_grad_strong_form}
	\int_{\holdalldomain}\left<(\operatorname{id} - A \Delta)^s \, \nabla_R L(\remainderdomain), \vW\right> \d \bm{x}  = \partial_{\remainderdomain} L(\remainderdomain,\state,\adjstate) [\vW] \quad \forall\, \vW\in C_0^\infty(\holdalldomain,\R^2).
\end{align}
or, analogously in weak form,
\begin{align}
	\label{eq:shape_grad_weak_form}
	\int_{\holdalldomain}\nabla_R L(\remainderdomain)^\top \vW + A \cdot \nabla \left(\nabla_R L(\remainderdomain)\right) : \nabla \vW \d \bm{x}  = \partial_{\remainderdomain} L(\remainderdomain,\state,\adjstate) [\vW] \quad \forall\, \vW\in C_0^\infty(\holdalldomain,\R^2).
\end{align}

Additionally, we also need to enforce the fracture irreversibility. %
 In order to guarantee the irreversibility of the fracture, we need to ensure that $\physdomain \subset \physdomain_{i}$ at each time step $t_i$, $i = 1, \ldots, N$, cf.~\eqref{eqn:OptProblemFinal}. Due to the complexity of enforcing this condition, it is replaced by an even more strict but easier-to-implement condition that the physical domain is not allowed to grow (i.e., the domain enclosed by the fracture is not allowed to shrink) between load steps if displaced by a vector field $\bm{V}$ as in~\eqref{eqn:DomainUpdate}. For this, we choose a standard frictionless contact problem enforced by a penalty approach. First, we describe the energy minimization problem of $F(\vV)=\frac{1}{2} a(\vV, \vV) - \shapeDeriv L(\remainderdomain, \state, \adjstate)[\vV]$ (i.e., a minimization problem different from and not related to~\eqref{eqn:OptProblemFinal}), where $a(\vV, \vV)$ is given as
 \begin{align}
 	\label{eq:EnergyFormulationRiemannianMetric}
 	a(\vV, \vV) = \int_{\holdalldomain} \vV^\top \vV + A \, {\nabla \vV} : {\nabla \vV} \d \bm{x},
 \end{align}
Assuming sufficient regularity, its first-order optimality criterion then yields equation~\eqref{eq:shape_grad_weak_form} for $s=1$. With $\bm{n}$ being the unit outward normal vector of $\physdomain$ it is immediately clear that a vector field~$\vV$, which fulfills $\vV^\top \bm{n} \leq -\epsilon$ on $\shape$, $\epsilon \geq 0$, complies with the condition of not allowing shrinkage of the physical domain. Therefore, the minimization problem to obtain the deformation equation reads
\begin{align*}
	\min_{\vV \in C_0^\infty(\holdalldomain,\R^2)} F(\vV) \quad \text{subject to} \quad \vV^\top \bm{n} \leq -\epsilon \text{ on } \shape.
\end{align*}
This can be reformulated to
\begin{align*}
	\min_{\vV \in C_0^\infty(\holdalldomain,\R^2)} F(\vV) \quad \text{subject to} \quad \max\left(0,\vV^\top \bm{n} + \epsilon \right) = 0 \text{ on } \shape.
\end{align*}
The enforcement of the constraint $\max\left(0,\vV^\top \bm{n} + \epsilon \right) = 0$ on $\shape$, i.e., $\max\left(0,\vV^\top \bm{n} + \epsilon \right) = 0 \ \forall \bm{x} \in \shape$, is, obviously, challenging.  According to measure theory, replacing this constraint by
\begin{align*}
	\int_\shape \max\left(0,\vV^\top \bm{n} + \epsilon \right) \d \bm{s} = 0
\end{align*}
implies that the number of points $\bm{x} \in \shape$ where $\vV^\top \bm{n} \leq -\epsilon$ is violated has a measure of zero, as a set with a measure of zero does not contribute to the value of the integral (cf., e.g.,~\cite{Elstrodt2018}). %
By using a penalty approach with penalty parameter $\psi > 0$ to enforce the constraint and with an exponent of $3$ to ensure twice-differentia\-bility, the minimization problem reads
\begin{align*}
	\min_{\vV \in C_0^\infty(\holdalldomain,\R^2)} \tunderbrace{F(\vV) +  \frac{\psi}{3} \int_\shape \max\left(0,\vV^\top \bm{n} + \epsilon \right)^3 \d \bm{s}}{F_\psi(\vV)}.
\end{align*}
A stationary point of this minimization problem has to fulfill the optimality condition
\begin{align}
	\begin{aligned}
	0 &= \partial_{\vV} F_\psi(\vV) [\vW] \\
	&= \partial_{\vV} F(\vV) [\vW] + \frac{\psi}{3} \int_\shape \partial_{\vV} \max\left(0,\vV^\top \bm{n} + \epsilon \right)^3 [\vW] \d \bm{s} \quad \forall \vW \in C_0^\infty(\holdalldomain,\R^2).
	\end{aligned}
	\label{eqn:OptimalityConditionDeformationEquation}
\end{align}
Due to the linearity of $a$ with respect to both arguments, the first term in~\eqref{eqn:OptimalityConditionDeformationEquation} reads
\begin{align*}
	\partial_{\vV} F(\vV) [\vW] &= \frac{1}{2} \partial_{\vV} a (\vV,\vV) [\vW] - \partial_{\vV} \left( \shapeDeriv L(\remainderdomain, \state, \adjstate)[\vV] \right) [\vW] \\
	&= a (\vV,\vW) - \shapeDeriv L(\remainderdomain, \state, \adjstate)[\vW] \\
	&=\int_{\holdalldomain} \vV^\top \vW + A \, {\nabla \vV} : {\nabla \vW} \d \bm{x} - \shapeDeriv L(\remainderdomain, \state, \adjstate)[\vW].
\end{align*}
Lastly, the derivative of the $\max$-term (the second term in~\eqref{eqn:OptimalityConditionDeformationEquation}) is given by
\begin{align*}
	\partial_{\vV} \max\left(0,\vV^\top \bm{n} + \epsilon \right)^3 [\vW] = 3 \max\left(0,\vV^\top \bm{n} + \epsilon \right)^2 \cdot \vW^\top \bm{n},
\end{align*}
therefore, the deformation equation to solve for $\vV \in C_0^\infty(\holdalldomain,\R^2)$ reads
\begin{align}
	\label{eq:DeformationEquationWithIrreversibility}
	\begin{aligned}
	0 = &\int_{\holdalldomain}\vV^\top \vW + A \, {\nabla \vV} : {\nabla \vW} \d \bm{x} - \shapeDeriv L(\remainderdomain, \state, \adjstate)[\vW] \\
	&\tunderbrace{+ \, \psi \int_\shape \max\left(0,\vV^\top \bm{n} + \epsilon \right)^2 \cdot \vW^\top \bm{n} \d \bm{s} \hspace*{1.7cm}}{R(\vV, \vW)}
	\end{aligned}
\end{align}
for all $\vW \in C_0^\infty(\holdalldomain,\R^2)$. %
This partial differential equation can be solved using Newton's method while increasing the penalty parameter~$\psi$. Therefore, we already determine the derivative of the residual~$R(\vV, \vW)$ with respect to~$\vV$ in direction~$\widetilde{\vW}$, which is given by
\begin{align*}
	&\partial_{\vV} R(\vV, \vW) [\widetilde{\vW}] \\
	&=  \int_{\holdalldomain} \widetilde{\vW}^\top \vW + A \, {\nabla \widetilde{\vW}} : {\nabla \vW} \d \bm{x} + 2\psi \int_\shape \max\left(0,\vV^\top \bm{n} + \epsilon \right) \cdot \vW^\top \bm{n} \cdot \widetilde{\vW}^\top \bm{n} \d \bm{s}.
\end{align*}

\section{Numerical investigations on a specific model}\label{sec:numerical_experiments}

In the following, we describe two common benchmark problems to evaluate the quality of fracture propagation simulations, cf.~\cite{Ambati2015,Mang2019}: the single-edge notched tension and the single-edge notched shear test. In \Cref{sec:modelFormulation} we first introduce the models and provide relevant material parameters for the simulation. \Cref{sec:ComputationRiemannianGradient} then describes the numerical implementation that is used to solve \eqref{eqn:OptProblemFinal}, and \Cref{sec:NumericalImplementation} then illustrates the discretizations of the domain, the investigated parameters and shows computational results for the tension test in \Cref{sec:NumericalExperimentsTensionTest} and the shear test in \Cref{sec:NumericalExperimentsShearTest}.

\subsection{Model formulation}\label{sec:modelFormulation}

\begin{figure}[tbp]
	\centering%
	\setlength\figureheight{.5\textwidth}%
	\setlength\figurewidth{.5\textwidth}%
		\newlength\svgwidth
		\setlength\svgwidth{\figurewidth}
		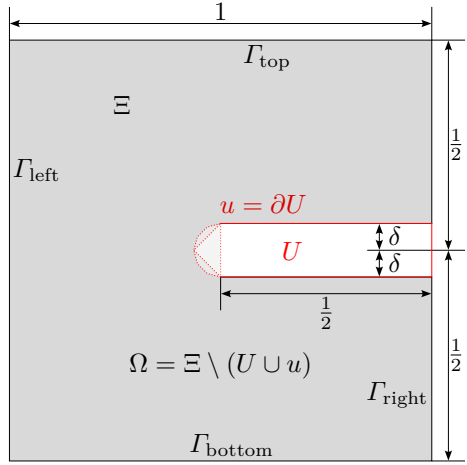
	\caption{Exact dimensions of the physical domain $\physdomain$ and the two-dimensional subdomain $\remainderdomain$ of the hold-all domain~$\holdalldomain$. The different types of fracture tips are also sketched.}%
	\label{fig:sketchModelDescription}%
\end{figure}

The hold-all domain $\holdalldomain \subset \R^2$ is split into a physical domain~$\physdomain$, the fracture~$\shape$, and the remainder~$\remainderdomain$, see \Cref{fig:sketchHoldAllDomainAndRemainderDomain}. The exact dimensions can be found in~\Cref{fig:sketchModelDescription}. As we want to investigate the influence of different initial fracture tips, we describe the tip of the fracture either as ``flat'', as ``pointy'' or as (an approximation of) ``round'', as also sketched in~\Cref{fig:sketchModelDescription}. The influence of the width of the initial fracture of $2\delta$ with $0 < \delta \ll \frac{1}{2}$ will also be investigated. %

The bulk material behavior is modeled as linear elasticity with homogeneous
material, for which the Lamé parameters read $\lambda =
\SI{121.15e3}{\newton \per \square \milli\meter}$ and $\mu =
\SI{80.77e3}{\newton \per \square \milli\meter}$. The boundary
conditions for linear elasticity are chosen as homogeneous Dirichlet on $\Gb$, homogeneous
Neumann on $\Gl \cup \Gr \cup \shape$, and an
inhomogeneous
Dirichlet boundary on $\Gt$ which depends on the specific
benchmark and increases in time (cf. \Cref{sec:FracturePropagationObjectiveFunctional}). As already elaborated for \eqref{eqn:OptProblemFinal_1} and similar to \cite{Suchan2023}, we estimate the value of the 1-dimensional
Hausdorff measure of the fracture~$\fracture$ as half the length of~$\shape$. In order to facilitate a comparison to alternative approaches,
we will evaluate the bulk energy $E_{\text{bulk}}$ and the boundary force $\bm{\tau}=(\tau_{x_1}, \tau_{x_2})^\top=\int_{\Gt} \stress(\state) \, \bm{n} \, \mathrm{d} s$ in the same direction as the applied inhomogeneous Dirichlet boundary condition $\state_D = \left( w_{x_1}, w_{x_2} \right)^\top$. Further, we compare the angle at which the fracture initially propagates from the notch, and the position of the fracture tip~$\bm{x}_{\text{split}}$ at an applied displacement of $\SI{1.5e-2}{\milli\meter}$ for the single-edge notched shear test.

The choice of $A$ in \eqref{eqn:ShapeGradientPushforward}, \eqref{eq:shape_grad_weak_form} and \eqref{eq:DeformationEquationWithIrreversibility} still has to be described. It is obvious that, the larger $A$ is chosen, the more the smoothing of the Poisson part in the deformation equation takes effect. In~\cite{LoayzaRomero2025}, values in the range of $A=0.01$ to $A=0.1$ are chosen. We choose $A=10$, which provided better results than smaller values for $A$ in our experiments. Additionally, to ensure that the fracture remains as initially chosen in the region where no fracture growth is expected, we further restrict the Riemannian gradient $\nabla_R L(\remainderdomain)$ to the space $\{ \vV \in C_0^\infty(\holdalldomain,\R^2) \colon \vV = \bm{0} \text{ on } \shape_{\text{horizontal}} \}$, where $\shape_{\text{horizontal}} \coloneqq \{ (x_1,x_2)^\top \in \shape \colon (x_2 = \frac{1}{2}-\delta \lor x_2 = \frac{1}{2} + \delta) \land x_1 \geq \frac{1}{2} \}$.

We are now ready to discuss the numerical approach to solve~\eqref{eqn:OptProblemFinal} in the following section. 

\subsection{Computation of the Riemannian gradient $\nabla_R L(\remainderdomain)$ and gradient descent algorithm} \label{sec:ComputationRiemannianGradient}

In order to find a solution to the minimization problem~\eqref{eqn:OptProblemFinal} for each time step $t_i$, $i=1,\ldots,N$, we would like to use the Riemannian gradient $\nabla_R L(\remainderdomain)$ to accelerate our progress towards the minimum (in comparison to gradient-free methods).
Therefore, we need to repeatedly compute a solution to the optimality system, i.e., the state equation~\eqref{eq:StateEquation}, the adjoint equation~\eqref{eq:AdjointEquation} and the design equation~\eqref{eqn:ShapeDerivative}. Thus, we first numerically solve the state equation (a linear PDE) to obtain the displacement~$\state$ on~$\physdomain$ using a classical Finite Element method with linear Lagrange elements and a direct solver (LU decomposition). In the next step, the adjoint equation (also a linear PDE) to obtain the adjoint displacement~$\adjstate$ is solved with the same method, element type and solver. In the last step, we would like to obtain a representative of our Riemannian gradient by solving~\eqref{eq:DeformationEquationWithIrreversibility}. We have to replace the space $C_0^\infty(\holdalldomain,\R^2)$ for $\nabla_R L(\remainderdomain)$ in \Cref{sec:modelFormulation} by a Sobolev space in order to use the Finite Element method to solve \eqref{eq:DeformationEquationWithIrreversibility} numerically. As we want to impose the least restrictions possible, the Sobolev space $\{ \vV \in H_0^1(\holdalldomain,\R^2) \colon \vV = \bm{0} \text{ on } \shape_{\text{horizontal}} \}$ is chosen because the (spatial) gradient of $\nabla_R L(\remainderdomain)$ needs to exist for \eqref{eq:shape_grad_weak_form} to be defined. Furthermore, computing the Riemannian gradient in~$\physdomain$ and on~$\shape$ is sufficient, as only the deformation of the boundary of~$\remainderdomain$ (which is also part of the boundary of~$\physdomain$ or does not change) affects any part of our optimization problem~\eqref{eqn:OptProblemFinal}. Now, we can numerically solve~\eqref{eq:DeformationEquationWithIrreversibility} using Newton's method with Armijo backtracking step size control to ensure a decrease in residual (initial step size $1$, backtracking parameter $10^{-4}$, step size halved if sufficient decrease condition not met). For the first iteration of the optimization of~\eqref{eqn:OptProblemFinal} we use an zero initial guess for Newton's method and increase the penalty factor by factors of $10$ from $\psi=10^{10}$ to $\psi=10^{15}$. For all other iterations we use the result of the previous optimization iteration as an initial guess for Newton's method and set $\psi=10^{15}$ immediately. In order to disincentivize too thin fractures which could cause the topology of the mesh to change, we set $\epsilon=10^{-7}$ in~\eqref{eq:DeformationEquationWithIrreversibility}.

Having obtained $\vV$ in $\physdomain$, we can now obtain a new domain using~\eqref{eqn:DomainUpdate}, i.e.,
\begin{align*}
	\physdomain_\tau = \{ \bm{x} \in \physdomain \colon \bm{x} + \tau \vV (\bm{x}) \}.
\end{align*}
In order to ensure a decrease in objective functional, the step length~$\tau$ in each optimization iteration is also obtained by Armijo backtracking step size control (initial step size $\SI{5e-3}{}$, backtracking parameter $10^{-4}$, step size halved if sufficient decrease condition not met). We assume convergence to a (local) minimum of the optimization in each time step when at least one of the three following conditions is met:
\begin{enumerate}
	\item The increase in fracture energy $E_{\text{frac}}$ (directly linked to the length increase of $\shape$) is below $10^{-8}$,
	\item the norm of $\vV$ that is induced by the outer Riemannian metric, i.e., $\sqrt{a(\vV,\vV)}$ with $a$ as in~\eqref{eq:EnergyFormulationRiemannianMetric}, falls below $10^{-4}$, or
	\item the Armijo backtracking step size control does not find a suitable step length~$\tau \geq 10^{-10}$.
\end{enumerate}
The whole simulation is stopped when the domain is (nearly) fully fractured or when a maximum number of loadsteps is reached.

With this knowledge, we can now focus on the software for generation of computational meshes and the numerical computation of the gradient.

\subsection{Numerical implementation} \label{sec:NumericalImplementation}

The numerical experiments are all performed using FEniCSx 0.8.0 \cite{Baratta2023}. The physical domain $\physdomain$ is discretized using Gmsh 4.13.1 \cite{Geuzaine2009} with an unstructured triangular mesh and without enforced mesh symmetry. Due to the large deformations of $\physdomain$, a remesher is required. If the mesh quality has deteriorated, an automatic remeshing procedure is used to generate a new discretization. The mesh quality is evaluated with the pyvista \cite{Sullivan2019} function \texttt{compute\_cell\_quality}. We use the minimum scaled Jacobian \cite{Knupp2000} over all triangular cells as an indication of the overall mesh quality. The quality threshold for remeshing is chosen as $30\%$ for all numerical experiments.

\begin{figure}[tbp]
	\centering%
	\setlength\figureheight{.25\textwidth}%
	\setlength\figurewidth{.25\textwidth}%
	\begin{subfigure}[t]{1.38\figurewidth}%
		\centering%
		\includetikz{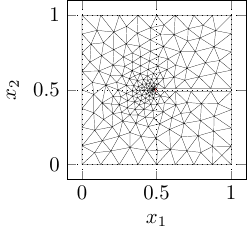}%
		\caption{Very coarse mesh containing 221 nodes and 394 triangular elements.}%
		\label{fig:VerycoarseMesh}%
	\end{subfigure}%
	\quad%
	\begin{subfigure}[t]{1.12\figurewidth}%
		\centering%
		\includetikz{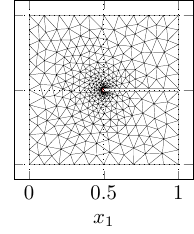}%
		\caption{Coarse mesh containing 403 nodes and 730 triangular elements.}%
		\label{fig:CoarseMesh}%
	\end{subfigure}%
	\quad%
	\begin{subfigure}[t]{1.12\figurewidth}%
		\centering%
		\includetikz{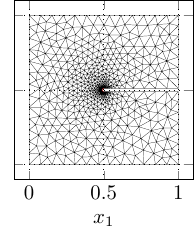}%
		\caption{Medium mesh containing 717 nodes and 1333 triangular elements.}%
		\label{fig:MediumMesh}%
	\end{subfigure}%
	\\[6pt]%
	\begin{subfigure}[t]{1.38\figurewidth}%
		\centering%
		\includetikz{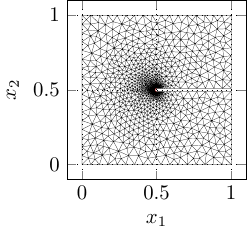}%
		\caption{Fine mesh containing 1322 nodes and 2513 triangular elements.}%
		\label{fig:FineMesh}%
	\end{subfigure}%
	\quad%
	\begin{subfigure}[t]{1.12\figurewidth}%
		\centering%
		\includetikz{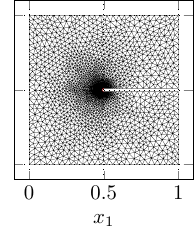}%
		\caption{Very fine mesh containing 3503 nodes and 6799 triangular elements.}%
		\label{fig:VeryfineMesh}%
	\end{subfigure}%
	\caption{Different refinement levels of the discretization of $\physdomain$ with the initially round fracture tip.}%
	\label{fig:MeshRefinementLevels}%
\end{figure}

\begin{figure}[tbp]
	\centering%
	\setlength\figureheight{.25\textwidth}%
	\setlength\figurewidth{.25\textwidth}%
	\begin{subfigure}[t]{1.3\figurewidth}%
		\centering%
		\includetikz{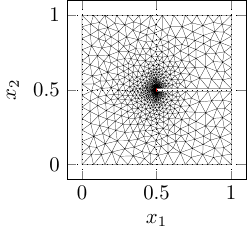}%
		\caption{Initial geometry with flat tip.}%
		\label{fig:MeshTipTestFlatTip}%
	\end{subfigure}%
	\quad%
	\begin{subfigure}[t]{\figurewidth}%
		\centering%
		\includetikz{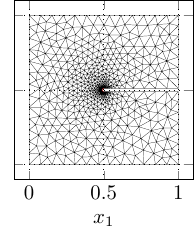}%
		\caption{Initial geometry with round tip.}%
		\label{fig:MeshTipTestRoundTip}%
	\end{subfigure}%
	\quad%
	\begin{subfigure}[t]{\figurewidth}%
		\centering%
		\includetikz{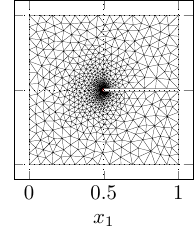}%
		\caption{Initial geometry with pointy tip.}%
		\label{fig:MeshTipTestPointyTip}%
	\end{subfigure}%
	\caption{Differently-formed initial notch tips with refinement level medium: flat, round or pointy.}%
	\label{fig:MeshTipTest}%
\end{figure}

\begin{figure}[tbp]
	\centering%
	\setlength\figureheight{.25\textwidth}%
	\setlength\figurewidth{.25\textwidth}%
	\begin{subfigure}[t]{1.3\figurewidth}%
		\centering%
		\includetikz{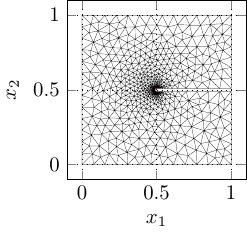}%
		\caption{Thicker fracture ($\delta=10^{-2}$).}%
		\label{fig:MeshFractureWidthTestThick}%
	\end{subfigure}%
	\quad%
	\begin{subfigure}[t]{\figurewidth}%
		\centering%
		\includetikz{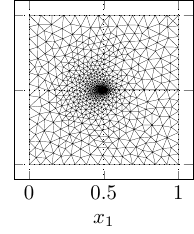}%
		\caption{Thinner fracture ($\delta=10^{-3}$).}%
		\label{fig:MeshFractureWidthTestThin}%
	\end{subfigure}%
	\caption{Different widths of the fracture: $\delta=10^{-2}$ and $\delta=10^{-3}$.}%
	\label{fig:MeshFractureWidthTest}%
\end{figure}

In order to review the influence of the mesh on the result, five different refinement levels of discretization have been generated: the very coarse mesh initially contains 221 nodes and 394 triangular elements; the coarse mesh initially contains 403 nodes and 730 triangular elements; the medium mesh initially contains 717 nodes and 1333 elements; the fine mesh initially contains 1322 nodes and 2513 elements; and the very fine mesh initially contains 3503 nodes and 6799 elements. For the round fracture tip and $\delta=10^{-2}$, the meshes are displayed in \Cref{fig:MeshRefinementLevels}.
In order to investigate the influence of the initial tip of the notch on the fracture propagation, three discretizations with different tips are available, which have already been described in \Cref{sec:modelFormulation} and sketched in~\Cref{fig:sketchModelDescription}. The discretizations for the medium mesh size can be seen in \Cref{fig:MeshTipTest}. Meshes with two different fracture widths, described by $\delta=10^{-2}$ and $\delta=10^{-3}$, and different values $\nu$ for the volume term, i.e., using $\nu=10$ and $\nu=100$, are also investigated. The different widths for the medium mesh size can be found in \Cref{fig:MeshFractureWidthTest}.

In the following, we will show different numerical results for the single-edge notched tension test in  \Cref{sec:NumericalExperimentsTensionTest} and the single-edge notched shear test in \Cref{sec:NumericalExperimentsShearTest}.

\subsubsection{Tension test} \label{sec:NumericalExperimentsTensionTest}

\begin{figure}[tbp]
	\centering%
	\setlength\figureheight{.49\textwidth}%
	\setlength\figurewidth{.49\textwidth}%
	\begin{subfigure}[t]{\figurewidth}%
		\centering%
		\newlength\svgwidth
		\setlength\svgwidth{\figurewidth}
		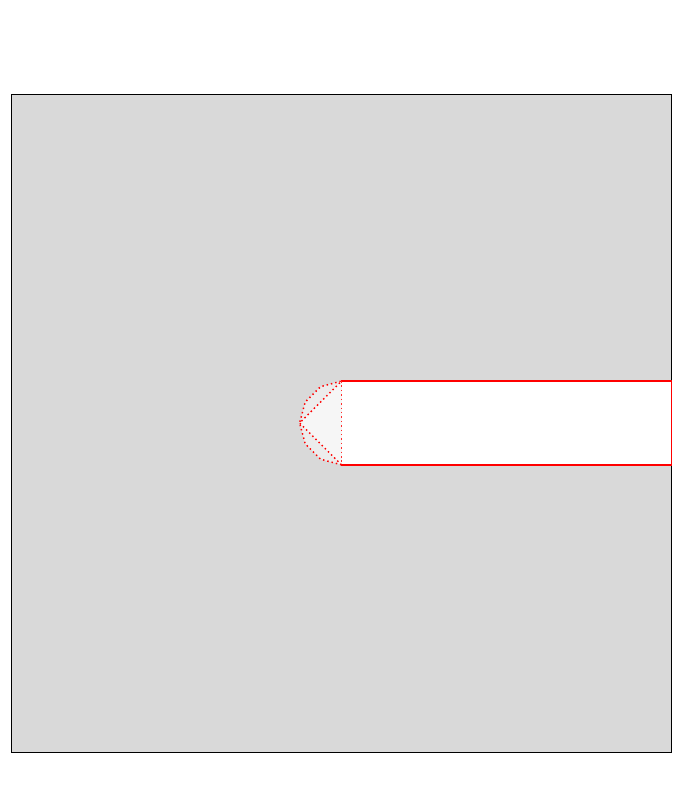
		\caption{Tension loading.}%
		\label{fig:sketchModelDescriptionLoadingTension}%
	\end{subfigure}%
	\hfill%
	\begin{subfigure}[t]{\figurewidth}%
		\centering%
		\newlength\svgwidth
		\setlength\svgwidth{\figurewidth}
		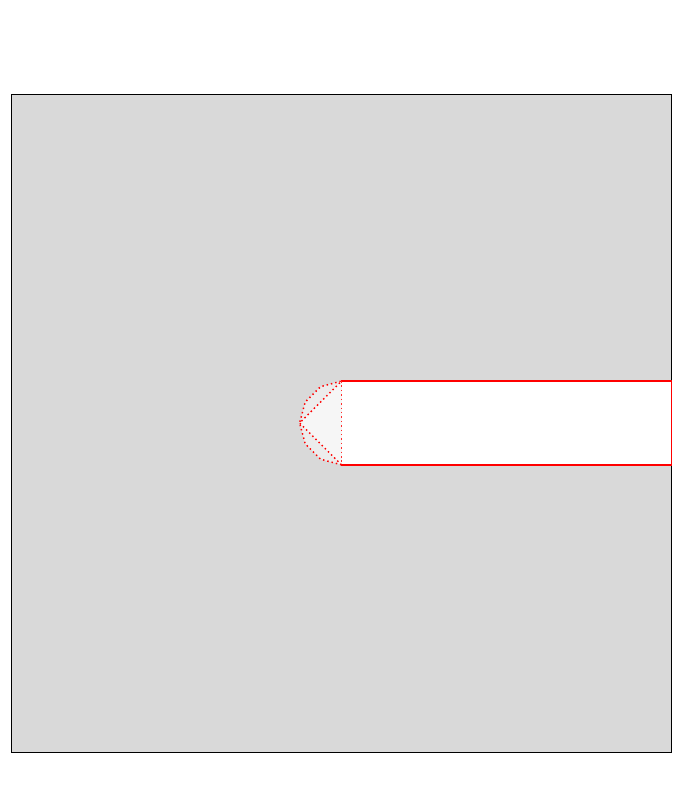
		\caption{Shear loading.}%
		\label{fig:sketchModelDescriptionLoadingShear}%
	\end{subfigure}%
	\caption{Physical domain $\physdomain$ with different types of applied loading: tension and shear.}%
	\label{fig:sketchModelDescriptionLoading}%
\end{figure}

For the numerical experiments in this section, we apply an inhomogeneous Dirichlet boundary condition on $\Gt$ in vertical direction as sketched in \Cref{fig:sketchModelDescriptionLoadingTension} which increases linearly in time (cf. \Cref{sec:FracturePropagationObjectiveFunctional}). The values are given in $\si{\micro\meter}$ unless specified otherwise. In order to reduce the numerical effort for this experiment at low displacements where no fracture propagation is expected, we use a displacement increment of $( 0, 0.5 )^\top$ until an inhomogeneous boundary condition of $(0, 3.5)^\top$ is reached, and then continue with a displacement increment of $( 0, 0.1 )^\top$. In a first step, we investigate the influence of mesh refinement. Then, the initial fracture tip shape is considered.  In a third numerical experiment, the effect of the fracture width, described by $\delta$, and the influence of the choice of $\nu$ on the fracture propagation result, are also studied.

\paragraph{Influence of the mesh refinement.} \label{sec:InfluenceMeshRefinementTension}

This numerical experiment investigates the influence of different mesh refinement levels on the fracture propagation. We compare simulations with the five different refinement levels in \Cref{fig:MeshRefinementLevels}. For this experiment, a round fracture tip is chosen, cf. \Cref{fig:sketchModelDescription} and \Cref{fig:MeshTipTestRoundTip}. Furthermore, an initial fracture width of $\delta=10^{-2}$ and a parameter of $\nu=10$ is selected.

\begin{figure}[tbp]
	\centering%
	\setlength\figureheight{.25\textwidth}%
	\setlength\figurewidth{.25\textwidth}%
	\begin{subfigure}[t]{1.38\figurewidth}%
		\centering%
		\includetikz{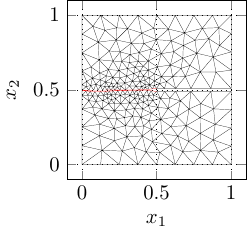}%
		\caption{Very coarse mesh.}%
		\label{fig:VerycoarseMeshFracturePathTension}%
	\end{subfigure}%
	\quad%
	\begin{subfigure}[t]{1.12\figurewidth}%
		\centering%
		\includetikz{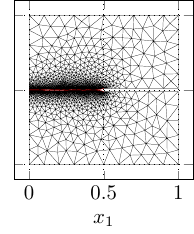}%
		\caption{Coarse mesh.}%
		\label{fig:CoarseMeshFracturePathTension}%
	\end{subfigure}%
	\quad%
	\begin{subfigure}[t]{1.12\figurewidth}%
		\centering%
		\includetikz{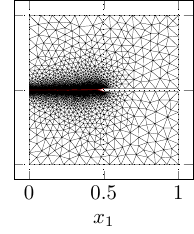}%
		\caption{Medium mesh.}%
		\label{fig:MediumMeshFracturePathTension}%
	\end{subfigure}%
	\\[6pt]%
	\begin{subfigure}[t]{1.38\figurewidth}%
		\centering%
		\includetikz{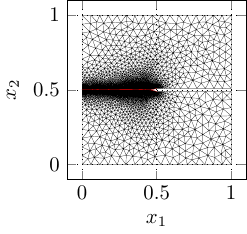}%
		\caption{Fine mesh.}%
		\label{fig:FineMeshFracturePathTension}%
	\end{subfigure}%
	\quad%
	\begin{subfigure}[t]{1.12\figurewidth}%
		\centering%
		\includetikz{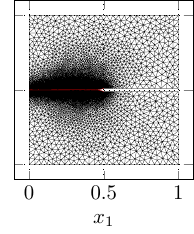}%
		\caption{Very fine mesh.}%
		\label{fig:VeryfineMeshFracturePathTension}%
	\end{subfigure}%
	\caption{Fracture paths for the tension test with different refinement levels of the discretization of $\physdomain$ with the initially round fracture tip.}%
	\label{fig:MeshRefinementLevelsFracturePathTension}%
\end{figure}

\begin{figure}[tbp]
	\centering%
	\setlength\figureheight{7cm}%
	\setlength\figurewidth{.48\textwidth}%
	\includelegend{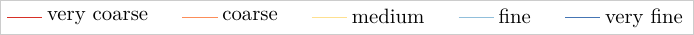}\\%
		\includetikz{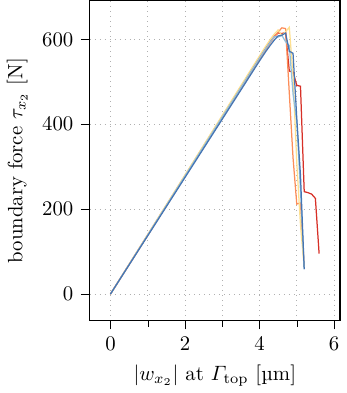}%
		\hfill%
		\includetikz{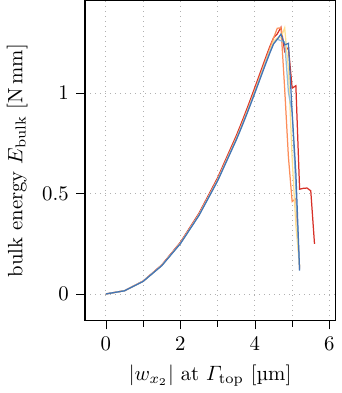}%
	\caption{Boundary force in $x_2$ direction and bulk energy at different displacements $w_{x_2}$ on $\Gt$ for the tension test with different refinement levels, a round fracture tip and $\delta=10^{-2}$.}%
	\label{fig:MeshRefinementTestTensionRoundForceAndEnergy}%
\end{figure}

The resulting fracture paths can be found in~\Cref{fig:MeshRefinementLevelsFracturePathTension}. All refinement levels provide perfectly horizontal fracture paths. When considering the boundary force in $x_2$ direction and energy stored in the bulk material in~\Cref{fig:MeshRefinementTestTensionRoundForceAndEnergy}, we observe very similar behavior over all refinement levels, except for the very coarse mesh, which shows slower mesh propagation than the other refinement levels. All refinement levels show the first major reduction in boundary force in the range at $\SI{4.8}{\micro\meter}$ or $\SI{4.9}{\micro\meter}$. The fracture has fully propagated through the domain at $\SI{5.6}{\micro\meter}$ for the very coarse mesh and $\SI{5.2}{\micro\meter}$ for the other meshes. As this type of fracture propagation is driven by stress, we associate the delay to the inaccuracy of stress computation, specifically at the fracture tip. In the other cases, we conclude very little influence of the mesh refinement.

\paragraph{Influence of the fracture tip.} \label{sec:InfluenceFractureTipTension}

\begin{figure}[tbp]
	\centering%
	\setlength\figureheight{.25\textwidth}%
	\setlength\figurewidth{.25\textwidth}%
	\begin{subfigure}[t]{1.3\figurewidth}%
		\centering%
		\includetikz{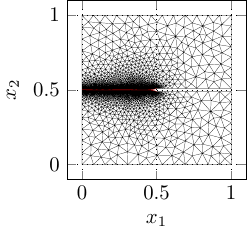}%
		\caption{Initially flat fracture tip.}%
		\label{fig:MeshTipTestTensionFlatTip}%
	\end{subfigure}%
	\quad%
	\begin{subfigure}[t]{\figurewidth}%
		\centering%
		\includetikz{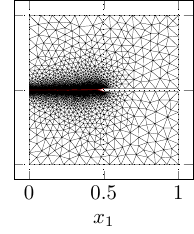}%
		\caption{Initially round fracture tip.}%
		\label{fig:MeshTipTestTensionRoundTip}%
	\end{subfigure}%
	\quad%
	\begin{subfigure}[t]{\figurewidth}%
		\centering%
		\includetikz{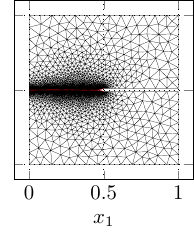}%
		\caption{Initially pointy fracture tip.}%
		\label{fig:MeshTipTestTensionPointyTip}%
	\end{subfigure}%
	\caption{Investigation of the fracture tip influence on the medium mesh (\Cref{fig:MediumMesh}) for the tension test.}%
	\label{fig:MeshTipTestTensionMesh}%
\end{figure}

\begin{figure}[tbp]
	\centering%
	\setlength\figureheight{7cm}%
	\setlength\figurewidth{.48\textwidth}%
	\includelegend{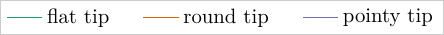}\\%
		\includetikz{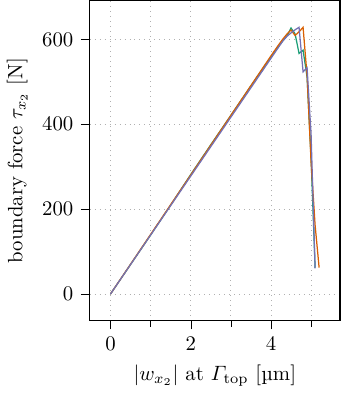}%
		\hfill%
		\includetikz{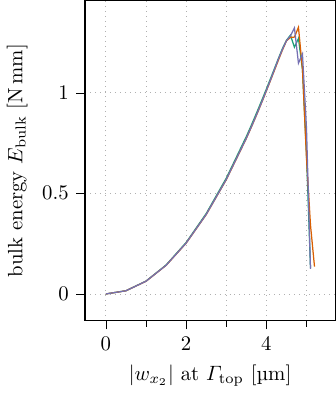}%
	\caption{Boundary force in $x_2$ direction and bulk energy at different displacements $w_{x_2}$ on $\Gt$ for the tension test with different initial tips.}%
	\label{fig:MeshTipTestTensionForceAndEnergy}%
\end{figure}

For the tension test, we investigate the influence of different
fracture tips on the fracture propagation simulation. We use the
medium mesh refinement level, $\delta=10^{-2}$ and
$\nu=10$. \Cref{fig:MeshTipTest} shows the initial configuration, to
which the increasing top boundary displacement is applied until the
domain is split in half. The fractures for the different initial tips
at the end of the simulation are shown in
\Cref{fig:MeshTipTestTensionMesh}. The simulations yield basically
identical fracture paths, independent of the initial tip shape. The
boundary force in $x_2$ direction and the bulk energy can be found in
\Cref{fig:MeshTipTestTensionForceAndEnergy}. At the initiation of
fracture propagation, small differences are visible in the boundary
forces at different applied displacements. However, the major
reduction in boundary force and bulk energy occurs at similar
displacements. %

\paragraph{Influence of the initial fracture width~$\delta$ and the volume parameter~$\nu$.} \label{sec:InfluenceFractureWidthAndVolumeParameterTension}

\begin{figure}[tbp]
	\centering%
	\setlength\figureheight{.25\textwidth}%
	\setlength\figurewidth{.25\textwidth}%
	\begin{subfigure}[t]{1.3\figurewidth}%
		\centering%
		\includetikz{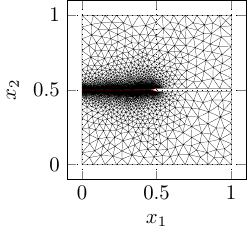}%
		\caption{Thicker fracture ($\delta=10^{-2}$), lower volume parameter ($\nu=10$).}%
		\label{fig:MeshFractureWidthAndVolumeParameterTestTensionThickNu10}%
	\end{subfigure}%
	\quad%
	\begin{subfigure}[t]{\figurewidth}%
		\centering%
		\includetikz{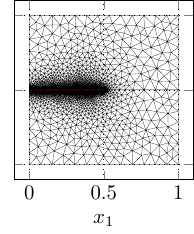}%
		\caption{Thinner fracture ($\delta=10^{-3}$), lower volume parameter ($\nu=10$).}%
		\label{fig:MeshFractureWidthAndVolumeParameterTestTensionThinWidthNu10}%
	\end{subfigure}%
	\quad%
	\begin{subfigure}[t]{\figurewidth}%
		\centering%
		\includetikz{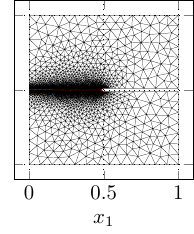}%
		\caption{Thicker fracture ($\delta=10^{-2}$), higher volume parameter ($\nu=100$).}%
		\label{fig:MeshFractureWidthAndVolumeParameterTestTensionThickWidthNu100}%
	\end{subfigure}%
	\caption{Investigation of the initial fracture width influence and the influence of the volume parameter on the medium mesh (\Cref{fig:MediumMesh}) for the tension test: $\delta=10^{-2}$ and $\delta=10^{-3}$; $\nu=10$ and $\nu=100$.}%
	\label{fig:MeshFractureWidthAndVolumeParameterTestTension}%
\end{figure}

\begin{figure}[tbp]
	\centering%
	\setlength\figureheight{7cm}%
	\setlength\figurewidth{.48\textwidth}%
	\includelegend{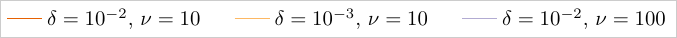}\\%
		\includetikz{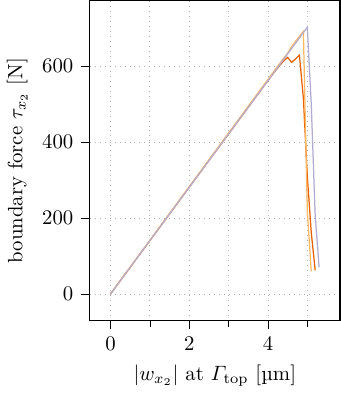}%
		\hfill%
		\includetikz{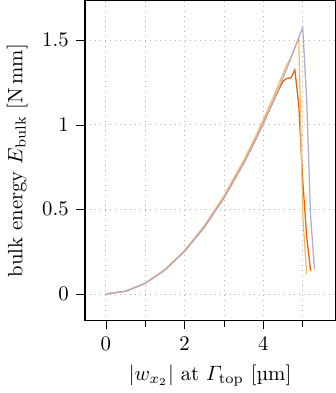}%
	\caption{Boundary force in $x_2$ direction and bulk energy at different displacements $w_{x_2}$ on $\Gt$ for the tension test with different initial fracture widths and volume parameter: $\delta=10^{-2}$ and $\delta=10^{-3}$; $\nu=10$ and $\nu=100$.}%
	\label{fig:MeshFractureWidthAndVolumeParameterTestTensionForceAndEnergy}%
\end{figure}

The third experiment for the single-edge notched tension test considers the influence of the initial fracture width $\delta$ and the volume parameter $\nu$. We choose the medium refinement level and a round fracture tip. An initial fracture width with $\delta=10^{-2}$ is compared to $\delta=10^{-3}$, and a volume parameter of $\nu=10$ is compared to $\nu=100$. The fracture paths are shown in~\Cref{fig:MeshFractureWidthAndVolumeParameterTestTension}. There is no notable difference in fracture path between $\delta=10^{-2}$ and $\delta=10^{-3}$ (\Cref{fig:MeshFractureWidthAndVolumeParameterTestTensionThickNu10} and \Cref{fig:MeshFractureWidthAndVolumeParameterTestTensionThinWidthNu10}), while a comparison between $\nu=10$ and $\nu=100$ (\Cref{fig:MeshFractureWidthAndVolumeParameterTestTensionThickNu10} and \Cref{fig:MeshFractureWidthAndVolumeParameterTestTensionThickWidthNu100}) shows a slight deviation from the purely horizontal fracture path with $\nu=100$ at closer inspection. The boundary force and bulk energy are displayed in~\Cref{fig:MeshFractureWidthAndVolumeParameterTestTensionForceAndEnergy}. We observe a later initiation of fracture propagation at $\SI{5.0}{\micro\meter}$ when using $\delta=10^{-3}$ in comparison to $\SI{4.9}{\micro\meter}$ for $\delta=10^{-2}$, while the fully fractured state is reached at a lower displacement ($\SI{5.1}{\micro\meter}$ vs. $\SI{5.2}{\micro\meter}$). Increasing the volume parameter requires a higher displacement of $\SI{5.1}{\micro\meter}$ for fracture initiation. A fully fractured domain is reached at $\SI{5.3}{\micro\meter}$. Thus, a thinner fracture is preferred, while a high volume parameter artificially delays the moment of initiation of fracture by disincentivizing the increase in domain volume. It should however be noted that a very thin fracture may only have very few elements at the fracture tip, which in turn can decrease the numerical accuracy of stress computations at the tip.

\subsubsection{Shear test} \label{sec:NumericalExperimentsShearTest}

We now continue with the analysis of the shape optimization approach with strain splitting for the shear test. This loading is sketched in \Cref{fig:sketchModelDescriptionLoadingShear}. The numerical experiments are again split into three parts. The same experiments are performed as for tension. We use a displacement increment of $( -0.5, 0 )^\top$ until an inhomogeneous boundary condition of $( -9, 0)^\top$ is reached, and then continue with a displacement increment of $( -0.1, 0 )^\top$. The maximum number of loadsteps is set to $100$, which yields a displacement of $( -17.2, 0 )^\top$ at the last loadstep. Furthermore, we compare the angle at which the fracture initially propagates from the fracture tip and the position of the fracture tip at a displacement of $( -15, 0)^\top$ to the results in~\cite{Ambati2015}.

\paragraph{Influence of the mesh refinement.} \label{sec:InfluenceMeshRefinementShear}

\begin{figure}[tbp]
	\centering%
	\setlength\figureheight{.25\textwidth}%
	\setlength\figurewidth{.25\textwidth}%
	\begin{subfigure}[t]{1.38\figurewidth}%
		\centering%
		\includetikz{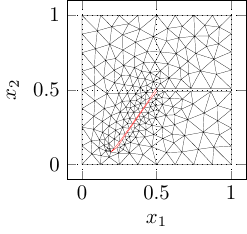}%
		\caption{Very coarse mesh.}%
		\label{fig:VerycoarseMeshFracturePathShear}%
	\end{subfigure}%
	\quad%
	\begin{subfigure}[t]{1.12\figurewidth}%
		\centering%
		\includetikz{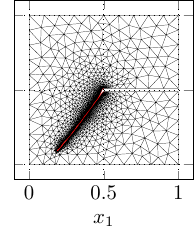}%
		\caption{Coarse mesh.}%
		\label{fig:CoarseMeshFracturePathShear}%
	\end{subfigure}%
	\quad%
	\begin{subfigure}[t]{1.12\figurewidth}%
		\centering%
		\includetikz{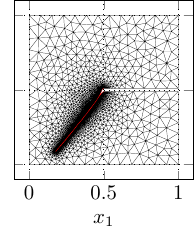}%
		\caption{Medium mesh.}%
		\label{fig:MediumMeshFracturePathShear}%
	\end{subfigure}%
	\\[6pt]%
	\begin{subfigure}[t]{1.38\figurewidth}%
		\centering%
		\includetikz{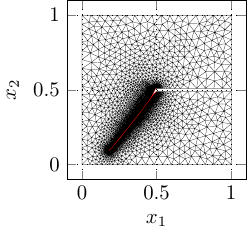}%
		\caption{Fine mesh.}%
		\label{fig:FineMeshFracturePathShear}%
	\end{subfigure}%
	\quad%
	\begin{subfigure}[t]{1.12\figurewidth}%
		\centering%
		\includetikz{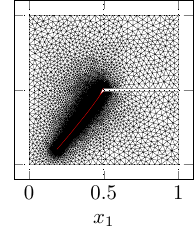}%
		\caption{Very fine mesh.}%
		\label{fig:VeryfineMeshFracturePathShear}%
	\end{subfigure}%
	\caption{Fracture paths for the shear test with different refinement levels of the discretization of $\physdomain$ with the initially round fracture tip after the last loadstep.}%
	\label{fig:MeshRefinementLevelsFracturePathShear}%
\end{figure}

\begin{figure}[tbp]
	\centering%
	\setlength\figureheight{7cm}%
	\setlength\figurewidth{.48\textwidth}%
	\includelegend{legend_mesh_refinement}\\%
		\includetikz{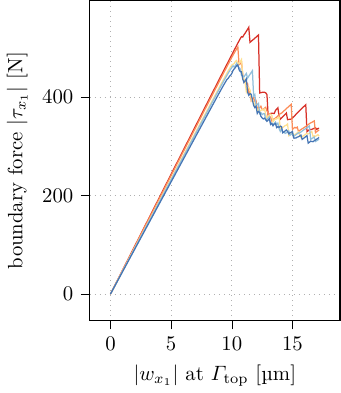}%
		\hfill%
		\includetikz{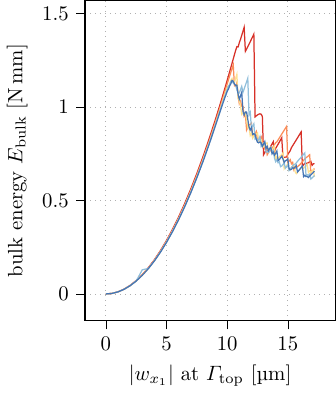}%
	\caption{Boundary force in $x_1$ direction and bulk energy at different displacements $w_{x_1}$ on $\Gt$ for the shear test with different refinement levels, a round fracture tip and $\delta=10^{-2}$.}%
	\label{fig:MeshRefinementTestShearRoundForceAndEnergy}%
\end{figure}

Similar to the tension test, we first compare the five different
refinement levels, which were already shown in
\Cref{fig:MeshRefinementLevels}. The initial width is chosen as
$\delta=10^{-2}$ and we have $\nu=10$ for the volume parameter. The
fractures after the last loadstep are shown
in~\Cref{fig:MeshRefinementLevelsFracturePathShear}. The fracture path
with all meshes is directed towards the bottom left.
In the plots for the boundary force and bulk energy, \Cref{fig:MeshRefinementTestShearRoundForceAndEnergy}, the displacements at fracture initiation are similar except for the very coarse mesh, which shows a delayed fracture initiation. In the later load steps the behavior is similar. Due to the delayed initiation of fracture we conclude that the very coarse mesh is not ideal for capturing the fracture propagation accurately, while the other meshes provide satisfactory results.

Furthermore, we would like to compare the angle at which the fracture
initially propagates and the position of the fracture tip for an
applied displacement at $\Gt$ of $\SI{-15}{\micro\meter}$ to
previous numerical results in \cite{Ambati2015}.
We list the values from \cite{Ambati2015} together with our results in \Cref{tab:ComparisonToAmbati}. It should be mentioned that the coordinate of the fracture tip was difficult to determine from \cite{Ambati2015} due to the nature of how phase-field methods describe the fracture. We have chosen to interpret the red area as fractured and tried to determine the position of the tip from the plot in the publication, knowing the size of the domain; nonetheless, this is only a rough approximation. A comparison of the results yields a similar initial angle of fracture propagation, with a difference of six to ten degrees. With a difference above $\SI{0.1}{\milli\meter}$, the coordinate of the fracture tip at a displacement of $\SI{-15}{\micro\meter}$ is not in accordance with the results in \cite{Ambati2015}. However, our numerical results do show very similar results when comparing them directly to each other, with coordinate differences of max. $\SI{0.023}{\milli\meter}$ in both $x_1$- and $x_2$-direction and angle differences of max. four degrees. We account the discrepancy between our results and \cite{Ambati2015} to the very different methods of computing stresses. We expect the stress in the phase field method to be lower than our computed stresses due to the dependence of the stress on the phase-field itself (cf.~\cite[Eq.~(43)]{Miehe2010}) when using the phase-field method. Since stresses are the cause of fracture propagation, we expect our fracture to start propagating at lower displacements, and therefore the fracture path to be longer at the same displacement. A similar effect can also be observed when comparing the tension test results from~\cite{Ambati2015} to our tension results in the previous section, i.e., \Cref{sec:NumericalExperimentsTensionTest}. This effect has already been observed in~\cite[Section~4.1]{Suchan2024}.

\begin{table}[tbp]
	\centering%
	\caption{Comparison to phase-field simulations in \cite{Ambati2015} of the $x_1$ and $x_2$ coordinate of the propagated fracture tip for the shear test at an applied displacement at $\Gt$ of $\SI{-15}{\micro\meter}$ and the angle of fracture growth shortly after initiation of fracture.}%
	\begin{tabular}{l|rrrrrr}%
		& \multicolumn{1}{c}{Ambati \cite{Ambati2015}} & \multicolumn{1}{c}{very coarse} & \multicolumn{1}{c}{coarse} & \multicolumn{1}{c}{medium} & \multicolumn{1}{c}{fine} & \multicolumn{1}{c}{very fine} \\
		\hline
		tip $x_1$ & 0.364 & 0.251 & 0.241 & 0.238 & 0.238 & 0.244 \\
		tip $x_2$ & 0.242 & 0.147 & 0.153 & 0.158 & 0.158 & 0.170 \\
		angle [$^\circ$] & 64  & 54  & 58  & 56  & 57  &  57
	\end{tabular}%
\label{tab:ComparisonToAmbati}%
\end{table}

\paragraph{Influence of the fracture tip.} \label{sec:InfluenceFractureTipShear}

\begin{figure}[tbp]
	\centering%
	\setlength\figureheight{.25\textwidth}%
	\setlength\figurewidth{.25\textwidth}%
	\begin{subfigure}[t]{1.3\figurewidth}%
		\centering%
		\includetikz{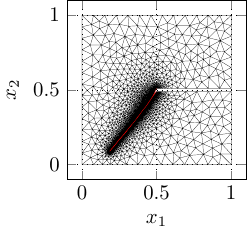}%
		\caption{Initially flat fracture tip.}%
		\label{fig:MeshTipTestShearFlatTip}%
	\end{subfigure}%
	\quad%
	\begin{subfigure}[t]{\figurewidth}%
		\centering%
		\includetikz{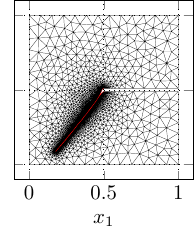}%
		\caption{Initially round fracture tip.}%
		\label{fig:MeshTipTestShearRoundTip}%
	\end{subfigure}%
	\quad%
	\begin{subfigure}[t]{\figurewidth}%
		\centering%
		\includetikz{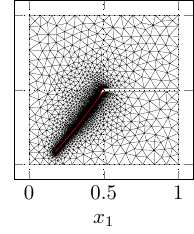}%
		\caption{Initially pointy fracture tip.}%
		\label{fig:MeshTipTestShearPointyTip}%
	\end{subfigure}%
	\caption{Investigation of the fracture tip influence on the medium mesh (\Cref{fig:MediumMesh}) for the shear test.}%
	\label{fig:MeshTipTestShearMesh}%
\end{figure}

\begin{figure}[tbp]
	\centering%
	\setlength\figureheight{7cm}%
	\setlength\figurewidth{.48\textwidth}%
	\includelegend{legend_fracture_tip}\\%
		\includetikz{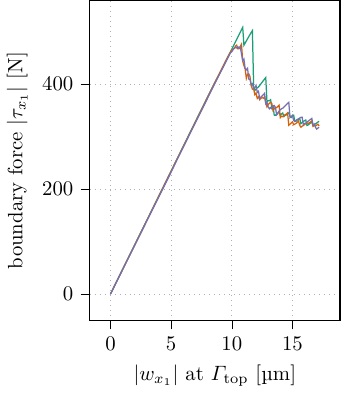}%
		\hfill%
		\includetikz{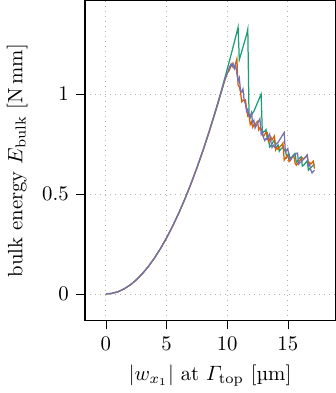}%
	\caption{Boundary force in $x_1$ direction and bulk energy at different displacements $w_{x_1}$ on $\Gt$ for the shear test with different initial tips.}%
	\label{fig:MeshTipTestShearForceAndEnergy}%
\end{figure}

For the shear test, we also investigate the different fracture tips (cf.~\Cref{fig:MeshTipTest}). We use the medium mesh refinement level, $\delta=10^{-2}$ and $\nu=10$. The obtained fracture paths are shown in~\Cref{fig:MeshTipTestShearMesh}. As in the tension test, no difference is visible for the fracture paths. The boundary force and bulk energy are depicted in~\Cref{fig:MeshTipTestShearForceAndEnergy}. The plots show the same qualitative behavior, with some differences specifically in the area of fracture initiation for the initially flat fracture tip. This is most likely related to the position of the sharp corners with respect to stress concentration in the shear test.

\paragraph{Influence of the initial fracture width~$\delta$ and the volume parameter~$\nu$.} \label{sec:InfluenceFractureWidthAndVolumeParameterShear}

\begin{figure}[tbp]
	\centering%
	\setlength\figureheight{.25\textwidth}%
	\setlength\figurewidth{.25\textwidth}%
	\begin{subfigure}[t]{1.3\figurewidth}%
		\centering%
		\includetikz{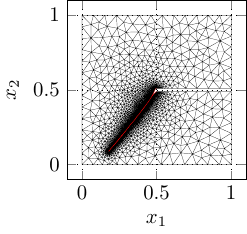}%
		\caption{Thicker fracture ($\delta=10^{-2}$), lower volume parameter ($\nu=10$).}%
		\label{fig:MeshFractureWidthAndVolumeParameterTestShearThickNu10}%
	\end{subfigure}%
	\quad%
	\begin{subfigure}[t]{\figurewidth}%
		\centering%
		\includetikz{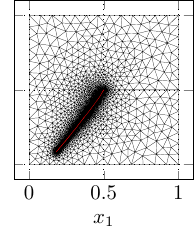}%
		\caption{Thin fracture ($\delta=10^{-3}$), lower volume parameter ($\nu=10$).}%
		\label{fig:MeshFractureWidthAndVolumeParameterTestShearThinWidthNu10}%
	\end{subfigure}%
	\quad%
	\begin{subfigure}[t]{\figurewidth}%
		\centering%
		\includetikz{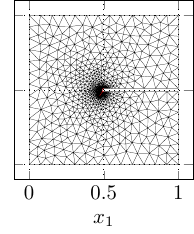}%
		\caption{Thicker fracture ($\delta=10^{-2}$), higher volume regularization ($\nu=100$).}%
		\label{fig:MeshFractureWidthAndVolumeParameterTestShearThickWidthNu100}%
	\end{subfigure}%
	\caption{Investigation of the initial fracture width influence and the influence of the volume parameter on the medium mesh (\Cref{fig:MediumMesh}) for the shear test: $\delta=10^{-2}$ and $\delta=10^{-3}$; $\nu=10$ and $\nu=100$.}%
	\label{fig:MeshFractureWidthAndVolumeParameterTestShear}%
\end{figure}

\begin{figure}[tbp]
	\centering%
	\setlength\figureheight{7cm}%
	\setlength\figurewidth{.48\textwidth}%
	\includelegend{legend_fracture_width_and_volume_reg}\\%
		\includetikz{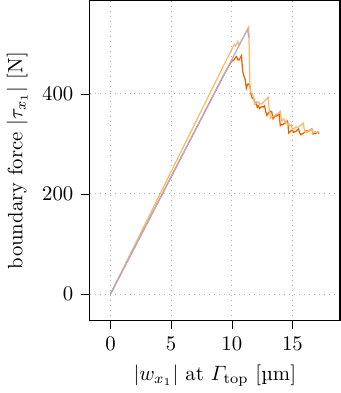}%
		\hfill%
		\includetikz{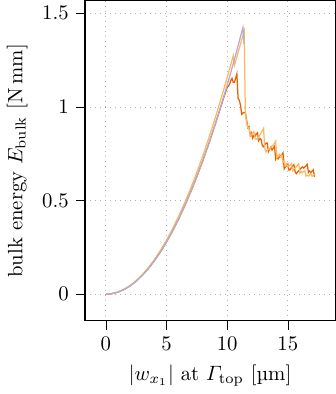}%
	\caption{Boundary force in $x_1$ direction and bulk energy at different displacements $w_{x_1}$ on $\Gt$ for the shear test with different initial fracture widths and volume parameter: $\delta=10^{-2}$ and $\delta=10^{-3}$; $\nu=10$ and $\nu=100$.}%
	\label{fig:MeshFractureWidthAndVolumeParameterTestShearForceAndEnergy}%
\end{figure}

The third numerical experiment with shear loading, as in the tension case, considers the influence of the initial fracture width $\delta$ and the volume parameter $\nu$. We again consider $\delta=10^{-2}$ or $\delta=10^{-3}$, and $\nu=10$ or $\nu=100$. The fracture paths for the different experiments can be found in~\Cref{fig:MeshFractureWidthAndVolumeParameterTestShear}. We first observe that the simulation for $\delta=10^{-2}$ and $\nu=100$ stopped prematurely due to a deteriorated mesh. Regarding the different fracture widths, we observe no notable difference in the fracture path. The boundary forces and bulk energy in~\Cref{fig:MeshFractureWidthAndVolumeParameterTestShearForceAndEnergy} show an approximately $5\%$ higher boundary force before initiation of fracture propagation for $\delta=10^{-3}$ in comparison to $\delta=10^{-2}$, most likely due to the additional bulk material that is present. The point of fracture initiation occurs at $\SI{-10.9}{\micro\meter}$ for $\delta=10^{-2}$ and $\SI{-11.5}{\micro\meter}$ for $\delta=10^{-3}$. After fracture initiation, the boundary force and bulk energy for $\delta=10^{-3}$ and $\delta=10^{-2}$ again look similar, with no clear advantage of one or the other thickness. Even though the experiment with $\nu=100$ failed prematurely, we can still determine that the fracture initiation is similar to $\delta=10^{-3}$ and $\nu=10$, which is similar to what was already observed for the tension test.

We have now described and evaluated all our numerical experiments. In the following section, we will conclude our paper and give a brief outlook for future investigations.

\section{Conclusion}\label{sec:summary_and_conclusion}

In this manuscript, we have used a previously-proposed approach to simulate fracture propagation by a shape optimization approach, which avoids the introduction of a length-scale parameter. We have provided an adaptation of the variational formulation of the fracture propagation problem to a shape optimization setting, while now also considering strain splitting to discard nonphysical fracture paths. This problem has been embedded into a Riemannian manifold framework by using a recently-proposed Riemannian shape space and an outer Riemannian metric, and the optimality system has been formulated. Furthermore, a novel approach to incorporate the irreversibility of the fracture has been described. Several numerical investigations for two classical benchmark problems, the single-edge notched tension and shear test have been performed, by which we have examined the influence of the mesh refinement, the initial fracture tip, the initial fracture width and the volume term. We have observed that very coarse discretizations do not provide satisfactory results, however coarser discretizations than in the phase field method can be used to generate satisfactory results. The influence of the refinement level (except for the very coarse mesh) has been shown to be low. In all cases, the horizontal direction of the fracture in the tension test has been observed. Furthermore, the initial fracture tip only has had a small influence on the fracture propagation. A thinner initial fracture can yield advantages, however it has to be ensured that the tip is still discretized with at least a few elements to maintain a sufficiently-accurate computation of the stresses at the fracture tip. A higher volume parameter has shown a delay in initiation of fracture, and the shear test has stopped prematurely due to a deteriorated mesh, and is therefore not recommended.
In summary, a shape optimization approach is a very good choice to simulate quasi-static, brittle fracture propagation problems, being robust to changes in parameters and levels of discretization refinement.

Further possible research directions include, but are not limited to, an extension to fractures in three-dimensional space, an alternative way of incorporating the irreversibility constraint of the optimization problem, and an adaptive strategy for the time discretization that accounts for the sudden initiation of fracture, especially in the tension test. Furthermore, alternative approaches to minimize the energy which models fracture propagation, e.g., a topology optimization approach, could be considered. It should also directly be possible to extend the optimization approach to time-dependent problems, e.g., cyclic loading or fatigue, however the failure criteria from cyclic loading or fatigue would have to be transferred to the shape optimization setting.

\section*{Acknowledgment}
The current work is a part of 
dtec.bw -- Digitalization and
Technology Research Center of the Bundeswehr within project ``Digitalisierung von Infrastrukturbauwerken zur Bauwerksüberwachung:
Structural Health Monitoring''. dtec.bw is funded by the European Union -- NextGenerationEU.

The authors are indebted to Lidiya Pryymak (Helmut Schmidt University/Uni\-ver\-sity of the Federal Armed Forces Hamburg, Germany) for many helpful comments and discussions
about the construction of the shape space.

\section*{Contribution}
\textbf{Tim Suchan:} Conceptualization, Methodology, Software, Formal analysis, Investigation, Writing -- Original Draft, Writing -- Review \& Editing, Visualization
\textbf{Winnifried Wollner:} Conceptualization, Methodology, Resources, Writing -- Original Draft, Writing -- Review \& Editing, Supervision
\textbf{Kathrin Welker:} Conceptualization, Methodology, Writing -- Original Draft, Writing -- Review \& Editing, Supervision, Funding acquisition

\crefalias{section}{appendix}
\appendix
\section{Numerical demonstration of strain splitting}\label{app:demonstration_strain_splitting}
\subsection{Tension test}

\begin{figure}[tbp]
	\setlength\figureheight{.4\textwidth}%
	\setlength\figurewidth{.4\textwidth}%
	\hspace*{3.2em}\includetikz{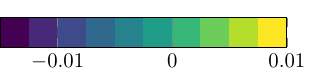}%
	\hspace*{0.1em}\includetikz{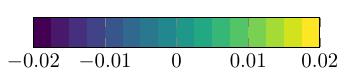}%
	\\%
	\includetikz{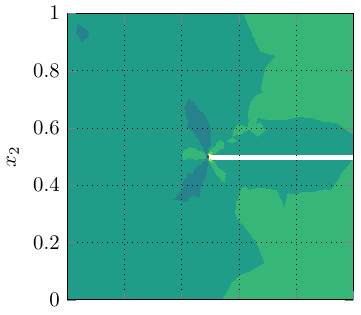}\qquad%
	\includetikz{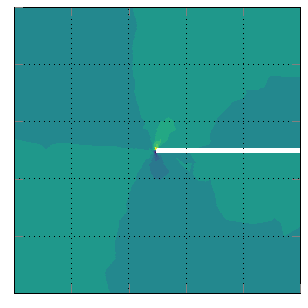}%
	\\%
	\hspace*{1.6em}\includetikz{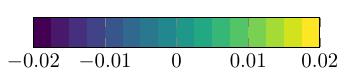}%
	\hspace*{2em}\includetikz{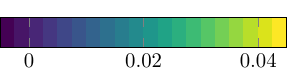}%
	\\%
	\includetikz{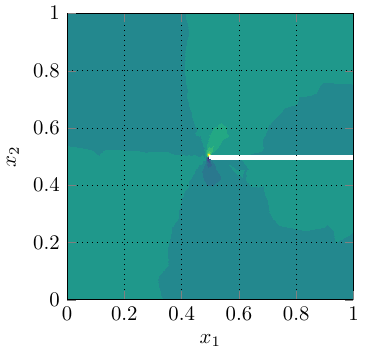}\qquad%
	\includetikz{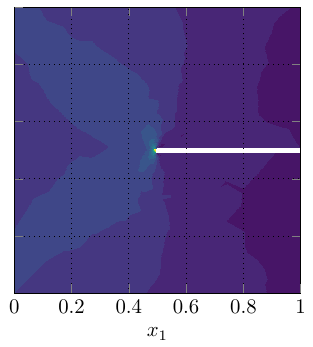}%
	\caption{Components of the strain tensor $\strain(\state)$ for the tension test with a top displacement of $\left(0, 5\right)^\top$. Top left: $\varepsilon_{1 1}(\state)$, top right: $\varepsilon_{1 2}(\state)$, bottom left: $\varepsilon_{2 1}(\state)$, bottom right: $\varepsilon_{2 2}(\state)$}%
	\label{fig:StrainTensorNumericalTension}%
\end{figure}

In order to numerically demonstrate this strain splitting method by finding the rotation angle $\alpha(\state)$, we consider the single-edge notched tension test (cf. \Cref{sec:NumericalExperimentsTensionTest}). We compute the displacement~$\state$ using linear Finite Elements with the medium mesh (cf.~\Cref{fig:MediumMesh}) while prescribing a top displacement of $\left(0, 5 \right)^\top \si{\micro\meter}$. Only in order to visualize the results, we then interpolate the results for $\strain(\state) = \frac{1}{2} \left( \nabla \state + {\nabla \state}^\top \right)$ back to a linear Finite Element space. In \Cref{fig:StrainTensorNumericalTension}, we show each individual component of $\strain(\state)$ over the whole domain, i.e., $\varepsilon_{1 1}(\state)$, $\varepsilon_{1 2}(\state)$, $\varepsilon_{2 1}(\state)$ and $\varepsilon_{2 2}(\state)$ (cf. \eqref{eq:strainTensorDefinition}).

\begin{figure}[tbp]
	\setlength\figureheight{.4\textwidth}%
	\setlength\figurewidth{.4\textwidth}%
	\hspace*{1.4em}\includetikz{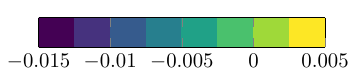}%
	\hspace*{-.5em}\includetikz{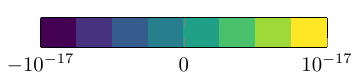}%
	\\%
	\includetikz{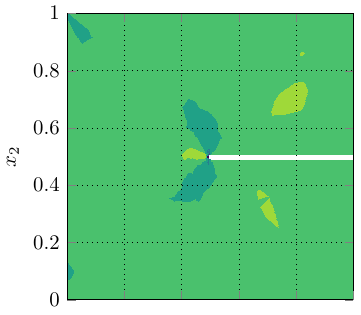}\qquad%
	\includetikz{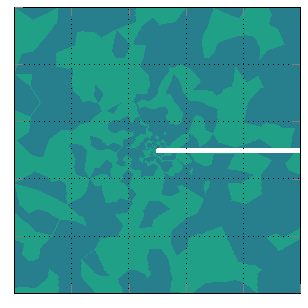}%
	\\%
	\hspace*{1.2em}\includetikz{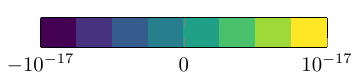}%
	\hspace*{1.2em}\includetikz{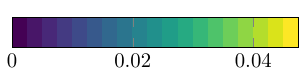}%
	\\%
	\includetikz{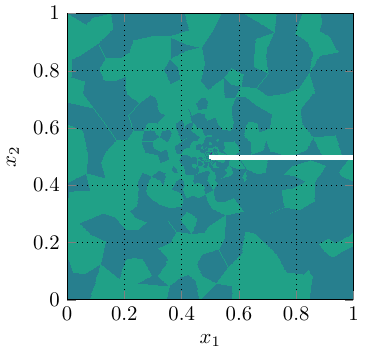}\qquad%
	\includetikz{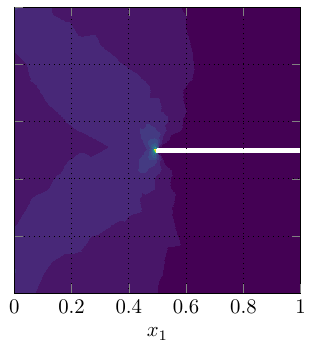}%
	\caption{Components of the rotated strain tensor $\bm{\Sigma}(\state)$ for the tension test with a top displacement of $\left(0, 5\right)^\top$. Top left: $\varepsilon_{1}(\state)$ (smaller eigenvalue of $\strain(\state)$), top right: $0$, bottom left: $0$, bottom right: $\varepsilon_{2}(\state)$ (larger eigenvalue of $\strain(\state)$)}
	\label{fig:StrainTensorRotatedNumericalTension}%
\end{figure}

In a first step, we check whether $\left| \varepsilon_{1 1}(\state) - \varepsilon_{2 2}(\state) \right| \leq 10^{-10}$ and $\left| \varepsilon_{1 2}(\state) \right| \leq 10^{-10}$. If this evaluates as true, then $\bm{Q}(\state)= \bm{I}$ and we already have a diagonal matrix. If this evaluates as false, then we compute the rotation angle $\alpha(\state)$ as described in \eqref{eq:RotationAngleAlpha}. This angle is then used in the rotation matrix $\bm{Q}(\state)$, cf. \eqref{eq:RotationMatrixQ}, to diagonalize the strain tensor as in \eqref{eq:StrainDecompositionDiagonalMatrix}, i.e., to obtain $\bm{\Sigma}(\state)$, where $\bm{Q}(\state)^{-1} = \bm{Q}(\state)^\top$, as $\bm{Q}(\state)$ is an orthogonal matrix by definition. For visualization, the rotated strain tensor is then again interpolated to a linear Finite Element space. \Cref{fig:StrainTensorRotatedNumericalTension} shows the numerical values of the entries in $\bm{\Sigma}(\state)$. From these plots, we conclude that $\bm{\Sigma}(\state)$ in fact is a diagonal matrix up to numerical accuracy. Furthermore, numerically evaluating $\varepsilon_{2}(\state)-\varepsilon_{1}(\state)$ indicates that $\varepsilon_{2}(\state)>\varepsilon_{1}(\state)$ in all of $\physdomain$.

\begin{figure}[tbp]
	\setlength\figureheight{.4\textwidth}%
	\setlength\figurewidth{.4\textwidth}%
	\hspace*{2.7em}\includetikz{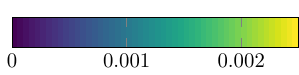}%
	\hspace*{2.6em}\includetikz{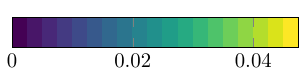}%
	\\%
	\includetikz{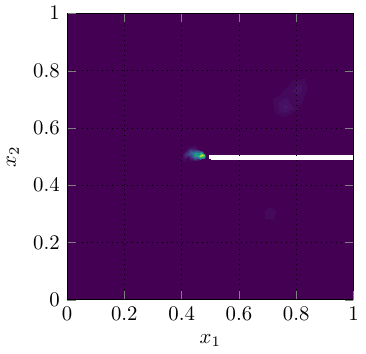}\qquad%
	\includetikz{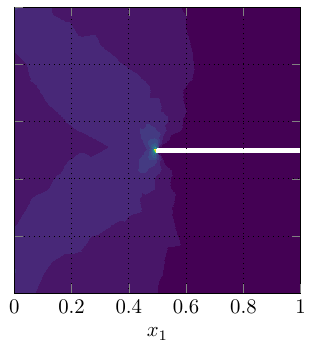}%
	\caption{Sorted nonnegative eigenvalues of the strain tensor $\strain(\state)$ for the tension test, obtained from the diagonal entries of $\bm{\Sigma}(\state)$. Left: $\max(0, \varepsilon_{1}(\state))$ (nonnegative smaller eigenvalue of $\strain(\state)$), right: $\max(0, \varepsilon_{2}(\state))$ (nonnegative larger eigenvalue of $\strain(\state)$)}%
	\label{fig:StrainTensorRotatedNonnegativeNumericalTension}%
\end{figure}

In a second step, we extract the diagonal entries from $\bm{\Sigma}(\state)$, i.e., the two eigenvalues $\varepsilon_{1}(\state)$ and $\varepsilon_{2}(\state)$, and numerically compute $\max(0, \varepsilon_{1}(\state))$ and $\max(0, \varepsilon_{2}(\state))$. This is again interpolated to linear Finite Elements for visualization and is shown in \Cref{fig:StrainTensorRotatedNonnegativeNumericalTension}. We observe numerical values above $-10^{-19}$, which is slightly negative, but within the range of numerical accuracy. This indicates that the nonnegative eigenvalues of $\strain(\state)$ for strain splitting have been successfully obtained with this method and can be used for later computations.

\subsection{Shear test}

\begin{figure}[tbp]
	\setlength\figureheight{.4\textwidth}%
	\setlength\figurewidth{.4\textwidth}%
	\hspace*{1.6em}\includetikz{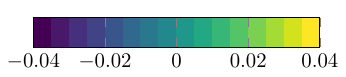}%
	\hspace*{1.8em}\includetikz{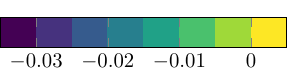}%
	\\%
	\includetikz{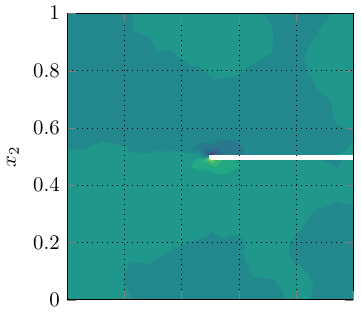}\qquad%
	\includetikz{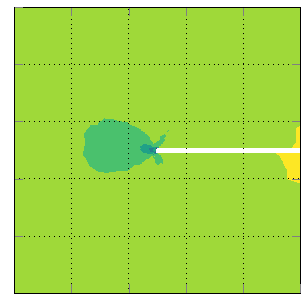}%
	\\%
	\hspace*{3.2em}\includetikz{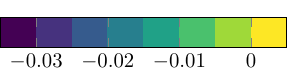}%
	\hspace*{3.2em}\includetikz{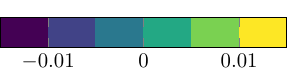}%
	\\%
	\includetikz{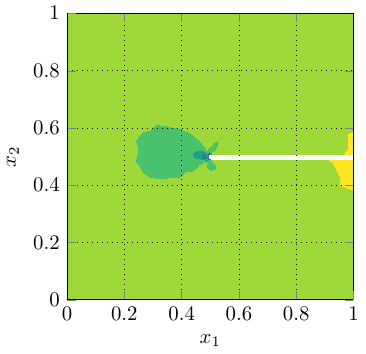}\qquad%
	\includetikz{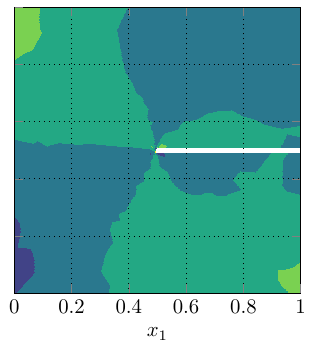}%
	\caption{Components of the strain tensor $\strain(\state)$ for the shear test with a top displacement of $\left(-10, 0\right)^\top$. Top left: $\varepsilon_{1 1}(\state)$, top right: $\varepsilon_{1 2}(\state)$, bottom left: $\varepsilon_{2 1}(\state)$, bottom right: $\varepsilon_{2 2}(\state)$}%
	\label{fig:StrainTensorNumericalShear}%
\end{figure}

We perform the same experiment for the single-edge notched shear test as for the tension test, with linear Finite Elements for the displacement~$\state$, medium refinement level, and a prescribed top displacement of $\left(-10, 0 \right)^\top \si{\micro\meter}$. \Cref{fig:StrainTensorNumericalShear} shows the components of the strain tensor, after interpolation to linear Finite Elements for visualization.

\begin{figure}[tbp]
	\setlength\figureheight{.4\textwidth}%
	\setlength\figurewidth{.4\textwidth}%
	\hspace*{3.2em}\includetikz{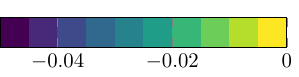}%
	\hspace*{0.45em}\includetikz{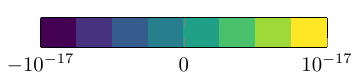}%
	\\%
	\includetikz{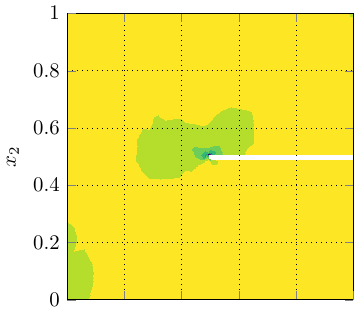}\qquad%
	\includetikz{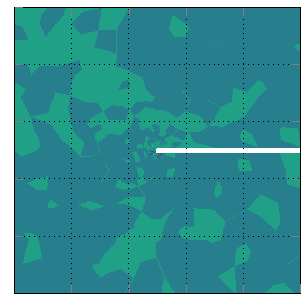}%
	\\%
	\hspace*{1.2em}\includetikz{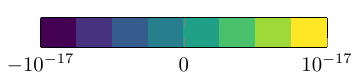}%
	\hspace*{1.2em}\includetikz{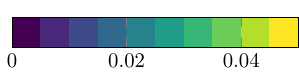}%
	\\%
	\includetikz{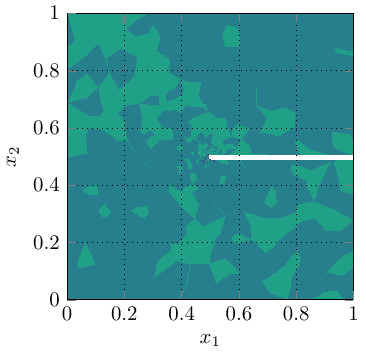}\qquad%
	\includetikz{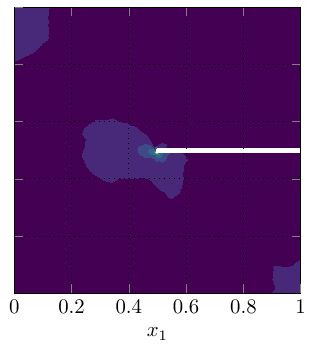}%
	\caption{Components of the rotated strain tensor $\bm{\Sigma}(\state)$ for the shear test with a top displacement of $\left(-10, 0\right)^\top$. Top left: $\varepsilon_{1}(\state)$ (smaller eigenvalue of $\strain(\state)$), top right: $0$, bottom left: $0$, bottom right: $\varepsilon_{2}(\state)$ (larger eigenvalue of $\strain(\state)$)}
	\label{fig:StrainTensorRotatedNumericalShear}%
\end{figure}

With the same procedure as for the tension test, we obtain the components of $\bm{\Sigma}(\state)$ and show them in \Cref{fig:StrainTensorRotatedNumericalShear}.

\begin{figure}[tbp]
	\setlength\figureheight{.4\textwidth}%
	\setlength\figurewidth{.4\textwidth}%
	\hspace*{1.3em}\includetikz{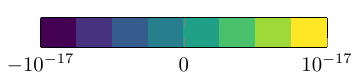}%
	\hspace*{1.1em}\includetikz{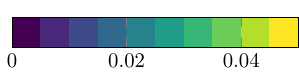}%
	\\%
	\includetikz{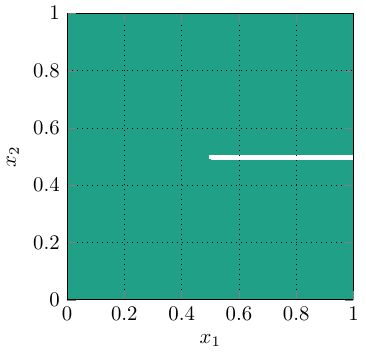}\qquad%
	\includetikz{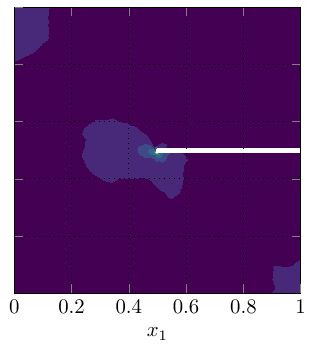}%
	\caption{Sorted nonnegative eigenvalues of the strain tensor $\strain(\state)$ for the shear test, obtained from the diagonal entries of $\bm{\Sigma}(\state)$. Left: $\max(0, \varepsilon_{1}(\state))$ (nonnegative smaller eigenvalue of $\strain(\state)$), right: $\max(0, \varepsilon_{2}(\state))$ (nonnegative larger eigenvalue of $\strain(\state)$)}%
	\label{fig:StrainTensorRotatedNonnegativeNumericalShear}%
\end{figure}

We again omit all negative eigenvalues, i.e., $\max(0, \varepsilon_{1}(\state))$ and $\max(0, \varepsilon_{2}(\state))$, and obtain \Cref{fig:StrainTensorRotatedNonnegativeNumericalShear}. Similar to the tension case, we again observe a very successful extraction of the nonnegative eigenvalues of $\strain(\state)$.

\section{Calculation of the state equation}\label{app:calculation_state_equation}

We observe the linear dependence of the state equation $c(\shape,\state,\adjstate)$ in~\eqref{eqn:OptProblemFinal_2} on the adjoint state~$\adjstate$. Therefore, the partial derivative of the Lagrange functional w.r.t. $\adjstate$ is given by
\begin{align*}
	\partial_{\adjstate} \left( c(\shape,\state(\shape),\adjstate(\shape)) \right) =
	\frac{\partial}{\partial \tau} \left. c(u, \state, \adjstate + \tau \tilde{\adjstate}) \right|_{\tau=0} = \frac{\partial}{\partial \tau} \left. \left( \int_\physdomain \stress(\state) : \strain(\adjstate + \tau \tilde{\adjstate}) \d \bm{x} \right) \right|_{\tau=0}.
\end{align*}
As $\strain(\adjstate + \tau \tilde{\adjstate}) = \strain(\adjstate) + \tau \strain(\tilde{\adjstate})$ this yields
\begin{align*}
	\frac{\partial}{\partial \tau} \left. c(u, \state, \adjstate + \tau \tilde{\adjstate}) \right|_{\tau=0} = \int_\physdomain \stress(\state) : \strain(\tilde{\adjstate}) \d \bm{x} \hastobe 0.
\end{align*}
Thus, the state equation is used to find a $\state \in H^1_D(\physdomain;\R^2)+ \state_D$ such that \eqref{eq:StateEquation} holds.

\section{Calculation of the adjoint equation}\label{app:calculation_adjoint_equation}

In order to obtain the adjoint equation we use the partial derivative of the Lagrange functional w.r.t. the state~$\state$, which is given as
\begin{align}
	\label{eqn:DerivativesForAdjointEquation}
	\begin{aligned}
		& \partial_{\state} \left( E_{\text{bulk}}(\shape,\state) + c(\shape,\state,\adjstate) \right) 
		= \frac{\partial}{\partial \tau} \left. \left( E_{\text{bulk}}(\shape,\state + \tau \tilde{\state}) + c(\shape,\state + \tau \tilde{\state},\adjstate) \right) \right|_{\tau=0} \\
		&= \frac{\partial}{\partial \tau} \left( \int_\physdomain \frac{\lambda}{2} \max(0,\trace(\strain(\state + \tau \tilde{\state})))^2 + \mu \trace \left( \bm{\Sigma}_{\text{max}}(\state + \tau \tilde{\state})^2 \right) \d \bm{x} \right. \\
		&\hphantom{=  \frac{\partial}{\partial \tau} \bigg(}\left. \left. + \int_\physdomain \stress(\state + \tau \tilde{\state}) : \strain(\adjstate) \d \bm{x} \right) \right|_{\tau=0}.
	\end{aligned}
\end{align}
As the piecewise derivatives of $\max(0,\cdot)^2$ coincide at $0$ here, we can now already give the derivative of the first term in \eqref{eqn:DerivativesForAdjointEquation}, which reads
\begin{align*}
	&\frac{\partial}{\partial \tau} \left. \left( \int_\physdomain \frac{\lambda}{2} \max(0,\trace(\strain(\state + \tau \tilde{\state})))^2 \d \bm{x} \right) \right|_{\tau=0} \\
	&= \int_\physdomain \frac{\lambda}{2} \cdot 2 \max(0,\trace(\strain(\state))) \cdot \frac{\partial}{\partial \tau} \left. \left( \trace(\strain(\state + \tau \tilde{\state})) \right) \right|_{\tau=0} \d \bm{x} \\
	&= \int_\physdomain \lambda \max(0,\trace(\strain(\state))) \cdot \trace(\strain(\tilde{\state})) \d \bm{x}.
\end{align*}
The third term is, due to the linearity of the stress tensor,
\begin{align*}
	\frac{\partial}{\partial \tau} \left. \left( \int_\physdomain \stress(\state + \tau \tilde{\state}) : \strain(\adjstate) \d \bm{x} \right) \right|_{\tau=0} = \int_\physdomain \stress(\tilde{\state}) : \strain(\adjstate) \d \bm{x} = \int_\physdomain \strain(\tilde{\state}) : \stress(\adjstate) \d \bm{x}.
\end{align*}
The second term requires some extra effort. First, we observe
\begin{align*}
	&\frac{\partial}{\partial \tau} \left. \left( \bm{\Sigma}_{\text{max}}(\state + \tau \tilde{\state})^2 \right) \right|_{\tau=0} \\
	&= \begin{pmatrix}
		\frac{\partial}{\partial \tau} \left. \max(0, \varepsilon_1(\state + \tau \tilde{\state}))^2 \right|_{\tau=0} & 0 \\
		0 & \frac{\partial}{\partial \tau} \left. \max(0, \varepsilon_2(\state + \tau \tilde{\state}))^2 \right|_{\tau=0}
	\end{pmatrix} \\
	&= \begin{pmatrix}
		2 \max(0, \varepsilon_1(\state)) \frac{\partial}{\partial \tau} \left. \varepsilon_1(\state + \tau \tilde{\state}) \right|_{\tau=0} & 0 \\
		0 & 2 \max(0, \varepsilon_2(\state)) \frac{\partial}{\partial \tau} \left. \varepsilon_2(\state + \tau \tilde{\state}) \right|_{\tau=0}
	\end{pmatrix} \\
	&= 2 \cdot \bm{\Sigma}_{\text{max}}(\state) \cdot \frac{\partial}{\partial \tau} \left.\bm{\Sigma}(\state + \tau \tilde{\state}) \right|_{\tau=0}
\end{align*}
We assume that the only possible nondifferentiability, that can occur when both eigenvalues coincide, can be neglected for our investigation. In the case $\varepsilon_{1 1}(\state) = \varepsilon_{2 2}(\state)$ and at the same time $\varepsilon_{1 2}(\state)=0$ (cf. \eqref{eqn:ConditionStrainDecomposition}), i.e., $\bm{\Sigma}(\state) = \strain(\state)$, we obtain
\begin{align*}
	\frac{\partial}{\partial \tau} \left. \left( \bm{\Sigma}_{\text{max}}(\state + \tau \tilde{\state})^2 \right) \right|_{\tau=0} = 2 \cdot \bm{\Sigma}_{\text{max}}(\state) \cdot \strain(\tilde{\state}),
\end{align*}
while in the opposite case, by using $\bm{\Sigma}(\state)=\bm{Q}(\state)^\top \, \strain(\state) \, \bm{Q}(\state)$ and the product rule, this gives
\begin{align*}
	&\frac{\partial}{\partial \tau} \left. \left( \bm{\Sigma}_{\text{max}}(\state + \tau \tilde{\state})^2 \right) \right|_{\tau=0} \\
	&= 2 \cdot \bm{\Sigma}_{\text{max}}(\state) \cdot \left( \frac{\partial}{\partial \tau} \left. \bm{Q}(\state + \tau \tilde{\state})^\top \right|_{\tau=0} \, \strain(\state) \, \bm{Q}(\state) + \bm{Q}(\state)^\top \, \strain(\tilde{\state}) \, \bm{Q}(\state) \right. \\
	& \hphantom{=2 \cdot \bm{\Sigma}_{\text{max}}(\state) \cdot \Big(\,} \left. + \, \bm{Q}(\state)^\top \, \strain(\state) \,  \frac{\partial}{\partial \tau} \left. \bm{Q}(\state + \tau \tilde{\state}) \right|_{\tau=0} \right).
\end{align*}
Then, using the chain rule, and with the derivative of $\alpha(\state)$ in \eqref{eq:RotationAngleAlpha} as
\begin{align*}
	\frac{\partial}{\partial \tau} \left. \alpha(\state + \tau \tilde{\state}) \right|_{\tau=0} &= \frac{\partial}{\partial \tau} \left. \frac{1}{2} \arctan\left(\xi(\state + \tau \tilde{\state})\right) \right|_{\tau=0} \\
	& \hphantom{= \ }+ \begin{cases}
		\frac{\partial}{\partial \tau} \left. 0 \right|_{\tau=0} & \text{, if } \varepsilon_{1 1}(\state) \leq \varepsilon_{2 2}(\state) \\
		\frac{\partial}{\partial \tau} \left. \frac{\pi}{2} \right|_{\tau=0} & \text{, else}
	\end{cases} \\
	&= \frac{\partial}{\partial \tau} \left. \frac{1}{2} \arctan\left(\xi(\state + \tau \tilde{\state})\right) \right|_{\tau=0},
\end{align*}
we obtain the derivative of $\bm{Q}(\state) = \bm{Q}(\alpha(\state))$ in \eqref{eq:RotationMatrixQ} as
\begin{align}
	\label{eqn:DerivationDerivativeQ}
	\begin{aligned}
		&\frac{\partial}{\partial \tau} \left. \bm{Q}(\state + \tau \tilde{\state}) \right|_{\tau=0} = \frac{\partial}{\partial \alpha(\state)} \bm{Q}(\alpha(\state)) \cdot  \frac{\partial}{\partial \tau} \left. \alpha(\state + \tau \tilde{\state}) \right|_{\tau=0} \\
		&= \frac{\partial}{\partial \alpha(\state)} \begin{pmatrix}
			\cos(\alpha(\state)) & -\sin(\alpha(\state)) \\
			\sin(\alpha(\state)) & \cos(\alpha(\state))
		\end{pmatrix} \cdot \frac{\partial}{\partial \tau} \left. \frac{1}{2} \arctan\left(\xi(\state + \tau \tilde{\state})\right) \right|_{\tau=0} \\
		&= \frac{\partial}{\partial \alpha(\state)} \begin{pmatrix}
			\cos(\alpha(\state)) & -\sin(\alpha(\state)) \\
			\sin(\alpha(\state)) & \cos(\alpha(\state))
		\end{pmatrix} \cdot \frac{\partial}{\partial \xi(\state)} \arctan\left(\xi(\state)\right) \frac{1}{2} \frac{\partial}{\partial \tau} \left. \xi(\state + \tau \tilde{\state}) \right|_{\tau=0} \\
		&= \tunderbrace{\begin{pmatrix}
				-\sin(\alpha(\state)) & -\cos(\alpha(\state)) \\
				\cos(\alpha(\state)) & -\sin(\alpha(\state))
		\end{pmatrix}}{\bm{R}(\state)} \cdot \frac{1}{1 + \xi(\state)^2} \cdot \frac{1}{2} \frac{\partial}{\partial \tau} \left. \xi(\state + \tau \tilde{\state}) \right|_{\tau=0}.
	\end{aligned}
\end{align}
Using the quotient rule and the linearity of the strain tensor w.r.t. $\state$, cf. \ref{app:calculation_state_equation}, we obtain
\begin{align*}
	\frac{\partial}{\partial \tau} \left. \xi(\state + \tau \tilde{\state}) \right|_{\tau=0} = \frac{2 \varepsilon_{1 2}(\tilde{\state})}{\varepsilon_{1 1}(\state)-\varepsilon_{2 2}(\state)} - \frac{2 \varepsilon_{1 2}(\state) \cdot \left(\varepsilon_{1 1}(\tilde{\state})-\varepsilon_{2 2}(\tilde{\state})\right)}{\left( \varepsilon_{1 1}(\state)-\varepsilon_{2 2}(\state) \right)^2},
\end{align*}
which yields the derivative of $\bm{Q}(\state)$ as
\begin{align*}
	&\frac{\partial}{\partial \tau} \left. \bm{Q}(\state + \tau \tilde{\state}) \right|_{\tau=0} \\
	&= \bm{R}(\state) \cdot \frac{1}{1 + \left(\frac{2 \varepsilon_{1 2}(\state)}{\varepsilon_{1 1}(\state)-\varepsilon_{2 2}(\state)}\right)^2} \cdot \left(\frac{\varepsilon_{1 2}(\tilde{\state})}{\varepsilon_{1 1}(\state)-\varepsilon_{2 2}(\state)} - \frac{\varepsilon_{1 2}(\state) \cdot \left(\varepsilon_{1 1}(\tilde{\state})-\varepsilon_{2 2}(\tilde{\state})\right)}{\left( \varepsilon_{1 1}(\state)-\varepsilon_{2 2}(\state) \right)^2}\right) \\
	&= \bm{R}(\state) \cdot \frac{\varepsilon_{1 2}(\tilde{\state}) \cdot \left( \varepsilon_{1 1}(\state)-\varepsilon_{2 2}(\state) \right) - \varepsilon_{1 2}(\state) \cdot \left(\varepsilon_{1 1}(\tilde{\state})-\varepsilon_{2 2}(\tilde{\state})\right)}{\left(\varepsilon_{1 1}(\state)-\varepsilon_{2 2}(\state)\right)^2 + 4 \varepsilon^2_{1 2}(\state)}.
\end{align*}
In summary, the derivative of the second term in \eqref{eqn:DerivativesForAdjointEquation} is given as
\begin{align*}
	&\frac{\partial}{\partial \tau} \left. \left( \int_\physdomain \mu \trace \left( \bm{\Sigma}_{\text{max}}(\state + \tau \tilde{\state})^2 \right) \d \bm{x} \right) \right|_{\tau=0} \\
	&= \int_\physdomain 2 \mu \begin{cases}
		\trace \left(\bm{\Sigma}_{\text{max}}(\state) \, \strain(\tilde{\state})\right) \d \bm{x} & \text{, if \eqref{eqn:ConditionStrainDecomposition}}, \\
			\frac{\varepsilon_{1 2}(\tilde{\state}) \cdot \left( \varepsilon_{1 1}(\state)-\varepsilon_{2 2}(\state) \right) - \varepsilon_{1 2}(\state) \cdot \left(\varepsilon_{1 1}(\tilde{\state})-\varepsilon_{2 2}(\tilde{\state})\right)}{\left(\varepsilon_{1 1}(\state)-\varepsilon_{2 2}(\state)\right)^2 + 4 \varepsilon^2_{1 2}(\state)} \\
			\ \ \cdot \trace \left( \bm{\Sigma}_{\text{max}}(\state) \left( \bm{R}(\state)^\top \strain(\state) \, \bm{Q}(\state) + \bm{Q}(\state)^\top \strain(\state) \, \bm{R}(\state) \right) \right) \\
			\ + \trace \left( \bm{\Sigma}_{\text{max}}(\state) \, \bm{Q}(\state)^\top \strain(\tilde{\state}) \, \bm{Q}(\state) \right) \d \bm{x}
		& \text{, else}
	\end{cases}
\end{align*}
and the adjoint equation follows from
$0 \hastobe \frac{\partial}{\partial \tau} \left. \left( E_{\text{bulk}}(\shape,\state + \tau \tilde{\state}) + h(\shape,\state + \tau \tilde{\state},\adjstate) \right) \right|_{\tau=0}$ %
as: Find $\adjstate \in H^1_D(\physdomain;\R^2)$ such that \eqref{eq:AdjointEquation} holds for all $\tilde{\state} \in H^1_D(\physdomain;\R^2)$.

\section{Calculation of the shape derivative}\label{app:calculation_shape_derivative}
We denote the material derivative of a function $p \colon \physdomain \times \R \to \R$ at a point $\bm{x} \in \physdomain$ in direction $\vW$ as $\deriv_m\left(p(\bm{x},0)\right)[\vW]$ or as $\dot{p}(\bm{x},0)[\vW]$, where the arguments are omitted if clear to avoid notational overhead. As described in~\cite{Berggren2009}, it is defined as
\begin{align}
	\label{eq:DefinitionMaterialDerivative}
	\deriv_m p(\bm{x},0) \left[\vW\right] \coloneqq \lim_{\tau \to 0^+} \frac{p(\bm{x}+\tau \vW, \tau) - p(\bm{x},0)}{\tau}.
\end{align}
Furthermore, the shape derivative of the integral of a function $p \colon \physdomain \times \R \to \R$ and its relation to the material derivative is also provided there, and is given as
\begin{align*}
	\partial_{\shape} \left( \int_\physdomain p(\bm{x},0) \d \bm{x} \right) \left[\vW\right] = \int_\physdomain \deriv_m p(\bm{x},0) \left[\vW\right] + \Div(\vW) \cdot p(\bm{x},0) \d \bm{x}.
\end{align*}
As described in \ref{app:calculation_adjoint_equation} we already have concluded differentiability of all terms in 
\begin{align*}
	E_{\text{bulk}}(\shape, \state) + E_{\text{frac}}(\shape) + E_{\text{reg}}(\shape) + c(\shape, \state, \adjstate),
\end{align*}
which extends to shape and material differentiability. We first calculate the material derivative of $\strain(\state)$ and show that the chain rule and the quotient rule for the material derivative hold, as it will be needed for the following material derivative calculations. With the material derivative of a gradient in \cite{Berggren2009} we have
\begin{align}
	\label{eq:MaterialDerivativeStrainTensor}
	\begin{aligned}
		&\deriv_m \left( \strain(\state) \right) = \frac{1}{2} \deriv_m \left(\nabla \state + {\nabla \state}^\top \right) = \frac{1}{2} \left( \left( \nabla \dot{\state} - \nabla \state \nabla \vW\right) + \left( \nabla \dot{\state} - \nabla \state \nabla \vW\right)^\top \right) \\
		&= \frac{1}{2} \left( \nabla \dot{\state} + {\nabla \dot{\state}}^\top \right) - \frac{1}{2} \left( \nabla \state \nabla \vW + \left(\nabla \state \nabla \vW\right)^\top \right) = \strain(\dot{\state}) - \sym\left( \nabla \state \nabla \vW \right).
	\end{aligned}
\end{align}
The material derivative of a composition $p=q \circ r$ of two differentiable functions $q \colon \physdomain \times \R \to \R$ and $r \colon \physdomain \times \R \to \physdomain \times \R$ is, by its definition in~\eqref{eq:DefinitionMaterialDerivative}, given by
\begin{align}
	\label{eqn:MaterialDerivativeChainRule}
	\begin{aligned}
		&\deriv_m p(\bm{x},0) \left[\vW\right] 
		= \lim_{\tau \to 0^+} \frac{p(\bm{x}+\tau \vW, \tau) - p(\bm{x},0)}{\tau} 
		= \lim_{\tau \to 0^+} \frac{q(r(\bm{x}+\tau \vW, \tau)) - q(r(\bm{x},0))}{\tau} \\
		&= \lim_{\tau \to 0^+} \frac{q(r(\bm{x}+\tau \vW, \tau)) - q(r(\bm{x},0))}{r(\bm{x}+\tau \vW, \tau) - r(\bm{x},0)} \cdot \frac{r(\bm{x}+\tau \vW, \tau) - r(\bm{x},0)}{\tau} \\
		&= \lim_{\tau \to 0^+} \frac{q(r(\bm{x}+\tau \vW, \tau)) - q(r(\bm{x},0))}{r(\bm{x}+\tau \vW, \tau) - r(\bm{x},0)} \cdot \lim_{\tau \to 0^+} \frac{r(\bm{x}+\tau \vW, \tau) - r(\bm{x},0)}{\tau} \\
		&= \lim_{r(\bm{x}+\tau \vW, \tau) \to r(\bm{x},0)} \frac{q(r(\bm{x}+\tau \vW, \tau)) - q(r(\bm{x},0))}{r(\bm{x}+\tau \vW, \tau) - r(\bm{x},0)} \cdot \lim_{\tau \to 0^+} \frac{r(\bm{x}+\tau \vW, \tau) - r(\bm{x},0)}{\tau} \\
		&= \partial_r q(r(\bm{x},0)) \cdot \deriv_m r(\bm{x},0) \left[\vW\right] .
	\end{aligned}
\end{align}
The material derivative of the function $p \colon \physdomain \times \R \to \R$, $p(\bm{x},\tau) = \frac{1}{q(\bm{x},\tau)}$ with \linebreak$q(\bm{x}+\tau \vW,\tau) \neq 0$ for all $\tau$ in a neighborhood of zero, can be derived as
\begin{align}
	\label{eq:MaterialDerivativeQuotientRule}
	\begin{aligned}
		&\deriv_m p(\bm{x},0) \left[\vW\right] = \lim_{\tau \to 0^+} \frac{p(\bm{x}+\tau \vW, \tau) - p(\bm{x},0)}{\tau} = \lim_{\tau \to 0^+} \frac{\frac{1}{q(\bm{x}+\tau \vW, \tau)} - \frac{1}{q(\bm{x},0)}}{\tau} \\
		&= \lim_{\tau \to 0^+} \frac{\frac{q(\bm{x},0)}{q(\bm{x},0) \cdot q(\bm{x}+\tau \vW, \tau)} - \frac{q(\bm{x}+\tau \vW, \tau)}{q(\bm{x},0) \cdot q(\bm{x}+\tau \vW, \tau)}}{\tau} = \lim_{\tau \to 0^+} \frac{q(\bm{x},0) - q(\bm{x}+\tau \vW, \tau)}{\tau \cdot q(\bm{x},0) \cdot q(\bm{x}+\tau \vW, \tau)} \\
		&= \lim_{\tau \to 0^+} \frac{-\left( q(\bm{x}+\tau \vW, \tau) - q(\bm{x},0) \right)}{\tau} \cdot \frac{1}{q(\bm{x},0) \cdot q(\bm{x}+\tau \vW, \tau)} \\
		&= \lim_{\tau \to 0^+} \frac{-\left( q(\bm{x}+\tau \vW, \tau) - q(\bm{x},0) \right)}{\tau} \cdot \lim_{\tau \to 0^+} \frac{1}{q(\bm{x},0) \cdot q(\bm{x}+\tau \vW, \tau)} \\
		&= - \deriv_m q(\bm{x},0) \left[\vW\right] \cdot \frac{1}{q(\bm{x},0)^2}.
	\end{aligned}
\end{align}
provided that $\deriv_m q(\bm{x},0) \left[\vW\right]$ exists and using that $\lim\limits_{\tau \to 0^+} q(\bm{x}+\tau \vW, \tau) = q(\bm{x},0)$.

We are now ready to calculate the shape derivative of \eqref{eq:LagrangeFunction}. It can be decomposed due to its linearity with respect to summation as
\begin{align}
	\label{eq:DesignEquationDecompositionSummation}
	\begin{aligned}
		&\shapeDeriv \left( L(\shape, \state, \adjstate) \right) \left[\vW\right] \\
		&= \shapeDeriv \left( \int_\physdomain \frac{\lambda}{2} \max(0,\trace(\strain(\state)))^2 \d \bm{x} \right) \left[\vW\right] + \shapeDeriv \left( \int_\physdomain \mu \trace \left( \bm{\Sigma}_{\text{max}}(\state)^2 \right) \d \bm{x} \right) \left[\vW\right] \\
		&\hphantom{=\ } + \shapeDeriv \left( \int_\shape \frac{G_c}{2} \d \bm{s} \right) \left[\vW\right] + \shapeDeriv \left( \int_{\remainderdomain} \nu \d \bm{x} \right) \left[\vW\right] + \shapeDeriv \left( \int_\physdomain \stress(\state) : \strain(\adjstate) \d \bm{x} \right) \left[\vW\right].
	\end{aligned}
\end{align}
We calculate the expression of the shape derivative with the material derivative approach from \cite{Berggren2009} and extensions in \cite{Geiersbach2023}.

\paragraph{First term in \eqref{eq:DesignEquationDecompositionSummation}.}
By observing that 
\begin{align*}
	\deriv_m \left( \trace(\strain(\state)) \right) &= \trace \left( \strain(\dot{\state}) - \sym\left( \nabla \state \nabla \vW \right) \right) = \trace \left( \strain(\dot{\state}) \right) - \trace \left( \sym\left( \nabla \state \nabla \vW \right) \right) \\
	&= \trace \left( \strain(\dot{\state}) \right) - \trace \left( \nabla \state \nabla \vW \right)
\end{align*}
we obtain
\begin{align*}
	&\shapeDeriv \left( \int_\physdomain \frac{\lambda}{2} \max(0,\trace(\strain(\state)))^2 \right) \left[\vW\right] \\
	&= \int_\physdomain \deriv_m \left( \frac{\lambda}{2} \max(0,\trace(\strain(\state)))^2 \right) + \Div(\vW) \cdot \frac{\lambda}{2} \max(0,\trace(\strain(\state)))^2 \d \bm{x} \\
	&= \int_\physdomain \frac{\lambda}{2} \cdot 2 \max(0,\trace(\strain(\state))) \, \deriv_m \left( \trace(\strain(\state)) \right) + \Div(\vW) \cdot \frac{\lambda}{2} \max(0,\trace(\strain(\state)))^2 \d \bm{x} \\
	&= \int_\physdomain \lambda \max(0,\trace(\strain(\state))) \cdot \trace \left( \strain(\dot{\state}) \right) - \lambda \max(0,\trace(\strain(\state))) \cdot \trace \left( \nabla \state \nabla \vW \right) \\
	&\hphantom{= \int_\physdomain \,}+ \Div(\vW) \cdot \frac{\lambda}{2} \max(0,\trace(\strain(\state)))^2 \d \bm{x}.
\end{align*}

\paragraph{Second term in \eqref{eq:DesignEquationDecompositionSummation}.}
The shape derivative of the second term is given as
\begin{align*}
	&\shapeDeriv \left( \int_\physdomain \mu \trace \left( \bm{\Sigma}_{\text{max}}(\state)^2 \right) \d \bm{x} \right) \left[\vW\right] \\
	&= \int_\physdomain \deriv_m \left( \mu \trace \left( \bm{\Sigma}_{\text{max}}(\state)^2 \right) \right) + \Div(\vW) \cdot \mu \trace \left( \bm{\Sigma}_{\text{max}}(\state)^2 \right) \d \bm{x}.
\end{align*}
With the material derivative of a second-order tensor as the material derivative of its entries, we can then write
\begin{align*}
	&\shapeDeriv \left( \int_\physdomain \mu \trace \left( \bm{\Sigma}_{\text{max}}(\state)^2 \right) \d \bm{x} \right) \left[\vW\right] \\
	&= \int_\physdomain  \mu \trace \left( \deriv_m \left( \bm{\Sigma}_{\text{max}}(\state)^2 \right) \right) + \Div(\vW) \cdot \mu \trace \left( \bm{\Sigma}_{\text{max}}(\state)^2 \right) \d \bm{x}.
\end{align*}
The material derivative of the tensor $\bm{\Sigma}_{\text{max}}(\state)^2$ is then described by
\begin{align*}
	\deriv_m \left( \bm{\Sigma}_{\text{max}}(\state)^2 \right) &= \begin{pmatrix}
		\deriv_m \left( \max(0, \varepsilon_1(\state))^2 \right) & 0 \\
		0 & \deriv_m \left( \max(0, \varepsilon_2(\state))^2 \right)
	\end{pmatrix} \\
	&= \begin{pmatrix}
		2 \max(0, \varepsilon_1(\state)) \cdot \deriv_m \left( \varepsilon_1(\state) \right) & 0 \\
		0 & 2 \max(0, \varepsilon_2(\state)) \cdot \deriv_m \left( \varepsilon_2(\state) \right)
	\end{pmatrix} \\
	&= 2 \cdot \bm{\Sigma}_{\text{max}}(\state) \cdot \deriv_m \left( \bm{\Sigma}(\state) \right).
\end{align*}
Similar to the derivation for the adjoint equation, we again observe two cases, which we will investigate individually:
\begin{enumerate}
	\item $\varepsilon_{1 1}(\state) = \varepsilon_{2 2}(\state) \land \varepsilon_{1 2}(\state)=0$ as in \eqref{eqn:ConditionStrainDecomposition} holds:\\[3pt]
	As $\bm{\Sigma}(\state) = \strain(\state)$, we immediately obtain 
	\begin{align*}
		\deriv_m \left( \bm{\Sigma}_{\text{max}}(\state)^2 \right) = 2 \cdot \bm{\Sigma}_{\text{max}}(\state) \cdot \left( \strain(\dot{\state}) - \sym\left( \nabla \state \nabla \vW \right) \right).
	\end{align*}
	\item $\varepsilon_{1 1}(\state) = \varepsilon_{2 2}(\state) \land \varepsilon_{1 2}(\state)=0$ does not hold:\\[3pt]
	We again use $\bm{\Sigma}(\state) = \bm{Q}(\state)^\top \, \strain(\state) \, \bm{Q}(\state)$ and therefore obtain
	\begin{align*}
		&\deriv_m \left( \bm{\Sigma}_{\text{max}}(\state)^2 \right) = 2 \cdot \bm{\Sigma}_{\text{max}}(\state) \cdot \deriv_m \left( \bm{Q}(\state)^\top \, \strain(\state) \, \bm{Q}(\state) \right) \\
		&= 2 \cdot \bm{\Sigma}_{\text{max}}(\state) \cdot \left( \deriv_m \left( \bm{Q}(\state)^\top \right) \, \strain(\state) \, \bm{Q}(\state) + \bm{Q}(\state)^\top \, \deriv_m \left(  \strain(\state) \right) \, \bm{Q}(\state) \right. \\
		&\hphantom{= 2 \cdot \bm{\Sigma}_{\text{max}}(\state) \cdot \big( \,}\left. + \, \bm{Q}(\state)^\top \, \strain(\state) \, \deriv_m \left( \bm{Q}(\state) \right) \right)
	\end{align*}
	using the product rule in~\cite{Berggren2009}. Using the chain rule for the material derivative of $\bm{Q}(\state)=\bm{Q}(\alpha(\state))=\bm{Q}(\alpha(\xi(\state)))$ in \eqref{eqn:MaterialDerivativeChainRule} we obtain
	\begin{align*}
		\deriv_m \left( \bm{Q}(\state) \right) = \tunderbrace{\partial_\alpha \bm{Q}(\alpha(\state))}{\bm{R}(\state)} \, \tunderbrace{\partial_\xi \alpha(\xi(\state))}{\frac{1}{2} \frac{1}{1 + \left(\frac{2 \varepsilon_{1 2}(\state)}{\varepsilon_{1 1}(\state)-\varepsilon_{2 2}(\state)}\right)^2}} \, \deriv_m \left( \xi(\state) \right).
	\end{align*}
	The two partial derivatives with respect to $\alpha$ and $\xi$ are the same as in~\eqref{eqn:DerivationDerivativeQ}. The material derivative of $\xi(\state)$, with the rule of the material derivative in \eqref{eq:MaterialDerivativeQuotientRule} and the product rule, is given by
	\begin{align*}
		\deriv_m \left( \xi(\state) \right) = \frac{2 \, \deriv_m \left( \varepsilon_{1 2}(\state) \right)}{ \varepsilon_{1 1}(\state) - \varepsilon_{2 2}(\state)} - \frac{2 \, \varepsilon_{1 2}(\state) \, \deriv_m \left( \varepsilon_{1 1}(\state) - \varepsilon_{2 2}(\state) \right)}{\left(\varepsilon_{1 1}(\state) - \varepsilon_{2 2}(\state)\right)^2}.
	\end{align*}
	From the material derivative of the strain tensor in \eqref{eq:MaterialDerivativeStrainTensor}, the material derivative of a partial derivative w.r.t. $x_i$, $i=1,2$ in~\cite{Berggren2009}, with $\state = \begin{pmatrix} w_1 , w_2 \end{pmatrix}^\top$ and $\vW = \begin{pmatrix} W_1 , W_2 \end{pmatrix}^\top$ we observe that
	\begin{align*}
		\deriv_m \left( \varepsilon_{1 1}(\state) \right) &= \deriv_m \left( \frac{\partial w_1}{\partial x_1} \right) = \frac{\partial \dot{w}_1}{\partial x_1} - \sum_{j=1}^2 \frac{\partial w_1}{\partial x_j} \frac{\partial W_j}{\partial x_1} = \varepsilon_{1 1}(\dot{\state}) - \sum_{j=1}^2 \frac{\partial w_1}{\partial x_j} \frac{\partial W_j}{\partial x_1} , \\
		\deriv_m \left( \varepsilon_{2 2}(\state) \right) &= \deriv_m \left( \frac{\partial w_2}{\partial x_2} \right) = \frac{\partial \dot{w}_2}{\partial x_2} - \sum_{j=1}^2 \frac{\partial w_2}{\partial x_j} \frac{\partial W_j}{\partial x_2} = \varepsilon_{2 2}(\dot{\state}) - \sum_{j=1}^2 \frac{\partial w_2}{\partial x_j} \frac{\partial W_j}{\partial x_2} , \\
		\deriv_m \left( \varepsilon_{1 2}(\state) \right) &= \frac{1}{2} \deriv_m \left( \frac{\partial w_1}{\partial x_2} + \frac{\partial w_2}{\partial x_1} \right) \\
		&= \frac{1}{2} \left( \frac{\partial \dot{w}_1}{\partial x_2} + \frac{\partial \dot{w}_2}{\partial x_1} \right)  - \frac{1}{2} \left(\sum_{j=1}^2 \frac{\partial w_1}{\partial x_j} \frac{\partial W_j}{\partial x_2} + \sum_{j=1}^2 \frac{\partial w_2}{\partial x_j} \frac{\partial W_j}{\partial x_1}\right) \\
		&= \varepsilon_{1 2}(\dot{\state}) - \frac{1}{2} \begin{pmatrix}
			\frac{\partial w_2}{\partial x_1} & \frac{\partial w_2}{\partial x_2} \\
			\frac{\partial w_1}{\partial x_1} & \frac{\partial w_1}{\partial x_2}
		\end{pmatrix} : \begin{pmatrix}
			\frac{\partial W_1}{\partial x_1} & \frac{\partial W_2}{\partial x_1} \\
			\frac{\partial W_1}{\partial x_2} & \frac{\partial W_2}{\partial x_2}
		\end{pmatrix} \\
		&= \varepsilon_{1 2}(\dot{\state}) - \frac{1}{2} \left( \begin{pmatrix}
			0 & 1 \\
			1 & 0
		\end{pmatrix} \, \begin{pmatrix}
			\frac{\partial w_1}{\partial x_1} & \frac{\partial w_1}{\partial x_2} \\
			\frac{\partial w_2}{\partial x_1} & \frac{\partial w_2}{\partial x_2}
		\end{pmatrix} \right) : \begin{pmatrix}
			\frac{\partial W_1}{\partial x_1} & \frac{\partial W_2}{\partial x_1} \\
			\frac{\partial W_1}{\partial x_2} & \frac{\partial W_2}{\partial x_2}
		\end{pmatrix} \\
		&= \varepsilon_{1 2}(\dot{\state}) - \frac{1}{2} \, \left( \begin{pmatrix}
			0 & 1 \\
			1 & 0
		\end{pmatrix} \, {\nabla \state}^\top \right) : \nabla \vW
	\end{align*}
	and we further obtain
	\begin{align*}
		&\deriv_m \left( \varepsilon_{1 1}(\state) - \varepsilon_{2 2}(\state) \right) = \varepsilon_{1 1}(\dot{\state}) - \varepsilon_{2 2}(\dot{\state}) - \sum_{j=1}^2 \frac{\partial w_1}{\partial x_j} \frac{\partial W_j}{\partial x_1} + \sum_{j=1}^2 \frac{\partial w_2}{\partial x_j} \frac{\partial W_j}{\partial x_2} \\
		&= \varepsilon_{1 1}(\dot{\state}) - \varepsilon_{2 2}(\dot{\state}) + \begin{pmatrix}
			-\frac{\partial w_1}{\partial x_1} & -\frac{\partial w_1}{\partial x_2} \\
			\frac{\partial w_2}{\partial x_1} & \frac{\partial w_2}{\partial x_2}
		\end{pmatrix} : \begin{pmatrix}
			\frac{\partial W_1}{\partial x_1} & \frac{\partial W_2}{\partial x_1} \\
			\frac{\partial W_1}{\partial x_2} & \frac{\partial W_2}{\partial x_2}
		\end{pmatrix} \\
		&= \varepsilon_{1 1}(\dot{\state}) - \varepsilon_{2 2}(\dot{\state}) + \left(  \begin{pmatrix}
			-1 & 0 \\
			0 & 1
		\end{pmatrix} \, \begin{pmatrix}
			\frac{\partial w_1}{\partial x_1} & \frac{\partial w_1}{\partial x_2} \\
			\frac{\partial w_2}{\partial x_1} & \frac{\partial w_2}{\partial x_2}
		\end{pmatrix} \right) : \begin{pmatrix}
			\frac{\partial W_1}{\partial x_1} & \frac{\partial W_2}{\partial x_1} \\
			\frac{\partial W_1}{\partial x_2} & \frac{\partial W_2}{\partial x_2}
		\end{pmatrix} \\
		&= \varepsilon_{1 1}(\dot{\state}) - \varepsilon_{2 2}(\dot{\state}) + \left( \begin{pmatrix}
			-1 & 0 \\
			0 & 1
		\end{pmatrix} \, {\nabla \state}^\top \right) : \nabla \vW .
	\end{align*}
	The material derivative in this case is therefore given as
	\begin{align*}
		&\deriv_m \left( \bm{\Sigma}_{\text{max}}(\state)^2 \right) \\
		&= 2 \cdot \bm{\Sigma}_{\text{max}}(\state) \cdot \bm{Q}(\state)^\top \left( \strain(\dot{\state}) - \sym\left( \nabla \state \nabla \vW \right) \right) \, \bm{Q}(\state) \\
		&\hphantom{=\ }+ 2 \cdot \frac{\left(\varepsilon_{1 2}(\dot{\state}) - \frac{1}{2} \, \left( \begin{pmatrix}
				0 & 1 \\
				1 & 0
			\end{pmatrix} \, {\nabla \state}^\top \right) : \nabla \vW \right) \cdot \left( \varepsilon_{1 1}(\state)-\varepsilon_{2 2}(\state) \right)}{\left(\varepsilon_{1 1}(\state)-\varepsilon_{2 2}(\state)\right)^2 + 4 \varepsilon^2_{1 2}(\state)} \\
		&\hphantom{= +\ \, }\cdot \bm{\Sigma}_{\text{max}}(\state) \cdot \left( \bm{R}(\state)^\top \strain(\state) \, \bm{Q}(\state) + \bm{Q}(\state)^\top \strain(\state) \, \bm{R}(\state) \right) \\
		&\hphantom{=\ }- 2 \cdot \frac{\varepsilon_{1 2}(\state) \cdot \left( \varepsilon_{1 1}(\dot{\state}) - \varepsilon_{2 2}(\dot{\state}) + \left( \begin{pmatrix}
				-1 & 0 \\
				0 & 1
			\end{pmatrix} \, {\nabla \state}^\top \right) : \nabla \vW \right)}{\left(\varepsilon_{1 1}(\state)-\varepsilon_{2 2}(\state)\right)^2 + 4 \varepsilon^2_{1 2}(\state)} \\
		&\hphantom{= +\ \, }\cdot \bm{\Sigma}_{\text{max}}(\state) \cdot \left( \bm{R}(\state)^\top \strain(\state) \, \bm{Q}(\state) + \bm{Q}(\state)^\top \strain(\state) \, \bm{R}(\state) \right).
	\end{align*}
\end{enumerate}
In summary, we obtain the material derivative of the second term as
\begin{align*}
	&\shapeDeriv \left( \int_\physdomain \mu \trace \left( \bm{\Sigma}_{\text{max}}(\state)^2 \right) \d \bm{x} \right) \left[\vW\right] \\
	&= \int_\physdomain \Div(\vW) \cdot \mu \trace \left( \bm{\Sigma}_{\text{max}}(\state)^2 \right) \\
	&\hphantom{=\int_\physdomain \ }+ 2 \mu \begin{cases}
		\trace \left( \bm{\Sigma}_{\text{max}}(\state) \cdot \left( \strain(\dot{\state}) - \sym\left( \nabla \state \nabla \vW \right) \right) \right) \d \bm{x} & \!\!\!\!\text{, if \eqref{eqn:ConditionStrainDecomposition}}, \\
			\frac{\left(\varepsilon_{1 2}(\dot{\state}) - \frac{1}{2} \, \left( \begin{pmatrix}
					0 & 1 \\
					1 & 0
				\end{pmatrix} \, {\nabla \state}^\top \right) : \nabla \vW \right) \cdot \left( \varepsilon_{1 1}(\state)-\varepsilon_{2 2}(\state) \right)}{\left(\varepsilon_{1 1}(\state)-\varepsilon_{2 2}(\state)\right)^2 + 4 \varepsilon^2_{1 2}(\state)} \\
			\ \ \cdot \trace \left( \bm{\Sigma}_{\text{max}}(\state) \cdot \left( \bm{R}(\state)^\top \strain(\state) \, \bm{Q}(\state) + \bm{Q}(\state)^\top \strain(\state) \, \bm{R}(\state) \right) \right) \\
			\ - \frac{\varepsilon_{1 2}(\state) \cdot \left( \varepsilon_{1 1}(\dot{\state}) - \varepsilon_{2 2}(\dot{\state}) + \left( \begin{pmatrix}
					-1 & 0 \\
					0 & 1
				\end{pmatrix} \, {\nabla \state}^\top \right) : \nabla \vW \right)}{\left(\varepsilon_{1 1}(\state)-\varepsilon_{2 2}(\state)\right)^2 + 4 \varepsilon^2_{1 2}(\state)} \\
			\ \ \cdot \trace \left( \bm{\Sigma}_{\text{max}}(\state) \cdot \left( \bm{R}(\state)^\top \strain(\state) \, \bm{Q}(\state) + \bm{Q}(\state)^\top \strain(\state) \, \bm{R}(\state) \right) \right) \\
			\ + \trace \left( \bm{\Sigma}_{\text{max}}(\state) \cdot \bm{Q}(\state)^\top \left( \strain(\dot{\state}) - \sym\left( \nabla \state \nabla \vW \right) \right) \, \bm{Q}(\state) \right) \d \bm{x}
		& \!\!\!\!\text{, else.}
	\end{cases}
\end{align*}
Ordering terms in this expression yields
\begin{align*}
	&\shapeDeriv \left( \int_\physdomain \mu \trace \left( \bm{\Sigma}_{\text{max}}(\state)^2 \right) \d \bm{x} \right) \left[\vW\right] \\
	&= \int_\physdomain \Div(\vW) \cdot \mu \trace \left( \bm{\Sigma}_{\text{max}}(\state)^2 \right) \\
	&\hphantom{=\int_\physdomain \ }+ 2 \mu \begin{cases}
		\trace \left( \bm{\Sigma}_{\text{max}}(\state) \cdot \strain(\dot{\state}) \right) - \trace \left( \bm{\Sigma}_{\text{max}}(\state) \cdot  \sym\left( \nabla \state \nabla \vW \right) \right) \d \bm{x} & \!\!\!\!\text{, if \eqref{eqn:ConditionStrainDecomposition}}, \\
			\trace \left( \bm{\Sigma}_{\text{max}}(\state) \cdot \bm{Q}(\state)^\top \strain(\dot{\state}) \, \bm{Q}(\state) \right) \\
			+\frac{\varepsilon_{1 2}(\dot{\state}) \cdot \left( \varepsilon_{1 1}(\state)-\varepsilon_{2 2}(\state) \right) - \varepsilon_{1 2}(\state) \cdot \left( \varepsilon_{1 1}(\dot{\state}) - \varepsilon_{2 2}(\dot{\state}) \right)}{\left(\varepsilon_{1 1}(\state)-\varepsilon_{2 2}(\state)\right)^2 + 4 \varepsilon^2_{1 2}(\state)} \\
			\ \ \cdot \trace \left( \bm{\Sigma}_{\text{max}}(\state) \cdot \left( \bm{R}(\state)^\top \strain(\state) \, \bm{Q}(\state) + \bm{Q}(\state)^\top \strain(\state) \, \bm{R}(\state) \right) \right) \\
			\ - \frac{ \left( \begin{pmatrix}
					-\varepsilon_{1 2}(\state) & \frac{\varepsilon_{1 1}(\state)-\varepsilon_{2 2}(\state)}{2} \\
					\frac{\varepsilon_{1 1}(\state)-\varepsilon_{2 2}(\state)}{2} & \varepsilon_{1 2}(\state)
				\end{pmatrix} {\nabla \state}^\top \right) : \nabla \vW}{\left(\varepsilon_{1 1}(\state)-\varepsilon_{2 2}(\state)\right)^2 + 4 \varepsilon^2_{1 2}(\state)} \\
			\ \ \cdot \trace \left( \bm{\Sigma}_{\text{max}}(\state) \cdot \left( \bm{R}(\state)^\top \strain(\state) \, \bm{Q}(\state) + \bm{Q}(\state)^\top \strain(\state) \, \bm{R}(\state) \right) \right) \\
			\ - \trace \left( \bm{\Sigma}_{\text{max}}(\state) \cdot \bm{Q}(\state)^\top \sym\left( \nabla \state \nabla \vW \right) \, \bm{Q}(\state) \right) \d \bm{x}
		& \!\!\!\!\text{, else.}
	\end{cases}
\end{align*}

\paragraph{Third term in \eqref{eq:DesignEquationDecompositionSummation}.}
Following~\cite[p.~80]{Sokolowski1992} the shape derivative of the third term directly follows as
\begin{align*}
	\shapeDeriv \left( \int_\shape \frac{G_c}{2} \d \bm{s} \right) \left[\vW\right] &= \frac{1}{2} \int_\shape {\nabla G_c}^\top \vW + G_c \left(\Div(\vW) - \bm{n}^\top \nabla \vW \bm{n}\right) \d \bm{s} \\
	&= \int_\shape \frac{G_c}{2} \left(\Div(\vW) - \bm{n}^\top \nabla \vW \bm{n}\right) \d \bm{s}.
\end{align*}

\paragraph{Forth term in \eqref{eq:DesignEquationDecompositionSummation}.}
With~\cite[p.~77]{Sokolowski1992}, the constant hold-all domain $\holdalldomain=\physdomain \cup \shape \cup \remainderdomain$ and $\nu$ being constant (specifically, with respect to $\bm{x}$) this shape derivative reads
\begin{align*}
	\shapeDeriv \left( \int_{\remainderdomain} \nu \d \bm{x} \right) \left[\vW\right] &= \shapeDeriv \left( \int_{\holdalldomain} \nu \d \bm{x} - \int_{\physdomain} \nu \d \bm{x} \right) \left[\vW\right] = \shapeDeriv \left( -\int_{\physdomain} \nu \d \bm{x} \right) \left[\vW\right] \\
	&= -\int_\physdomain \nu \Div(\vW) \d \bm{x}.%
\end{align*}

\paragraph{Fifth term in \eqref{eq:DesignEquationDecompositionSummation}.}
This shape derivative of the fifth term can be expressed with respect to the material derivative as
\begin{align*}
	&\shapeDeriv \left( \int_\physdomain \stress(\state) : \strain(\adjstate) \d \bm{x} \right) \left[\vW\right] = \int_\physdomain \deriv_m \left( \stress(\state) : \strain(\adjstate) \right) + \Div(\vW) \cdot \stress(\state) : \strain(\adjstate) \d \bm{x} \\
	&= \int_\physdomain \deriv_m \left( \stress(\state) \right) : \strain(\adjstate) + \stress(\state) : \deriv_m \left(\strain(\adjstate) \right) + \Div(\vW) \cdot \stress(\state) : \strain(\adjstate) \d \bm{x}
\end{align*}
and with $\stress(\state) = \mathbb{C} : \strain(\state)$ and therefore
\begin{align*}
	\deriv_m \left( \stress(\state) \right) &= \deriv_m \left( \mathbb{C} : \strain(\state) \right) = \mathbb{C} : \deriv_m \left( \strain(\state) \right)
	= \mathbb{C} : \left( \strain(\dot{\state}) - \sym\left( \nabla \state \nabla \vW \right) \right), %
\end{align*}
this yields
\begin{align*}
	&\shapeDeriv \left( \int_\physdomain \stress(\state) : \strain(\adjstate) \d \bm{x} \right) \left[\vW\right] \\
	&= \int_\physdomain \left( \mathbb{C} : \left( \strain(\dot{\state}) - \sym\left( \nabla \state \nabla \vW \right) \right) \right) : \strain(\adjstate) + \stress(\state) : \left( \strain(\dot{\adjstate}) - \sym\left( \nabla \adjstate \nabla \vW \right) \right) \\
	&\hphantom{= \int_\physdomain \,}+ \Div(\vW) \cdot \stress(\state) : \strain(\adjstate) \d \bm{x} \\
	&= \int_\physdomain \strain(\dot{\state}) : \stress(\adjstate) + \stress(\state) : \strain(\dot{\adjstate}) - \sym\left( \nabla \state \nabla \vW \right) : \stress(\adjstate) - \stress(\state) : \sym\left( \nabla \adjstate \nabla \vW \right) \\
	&\hphantom{= \int_\physdomain\,}+ \Div(\vW) \cdot \stress(\state) : \strain(\adjstate) \d \bm{x}.
\end{align*}

\paragraph{Shape derivative of the Lagrange functional}
Summation of all terms as described in \eqref{eq:DesignEquationDecompositionSummation} and sorting together the remaining material derivative terms gives 
\begin{align*}
	&\!\!\!\!\shapeDeriv \left( L(\shape, \state, \adjstate) \right) \left[\vW\right]\\
	=& \int_\physdomain \stress(\state) : \strain(\dot{\adjstate}) \d \bm{x} + \int_\physdomain \lambda \max(0,\trace(\strain(\state))) \cdot \trace \left( \strain(\dot{\state}) \right) +  \strain(\dot{\state}) : \stress(\adjstate) \d \bm{x} \\
	&+ \int_\physdomain 2 \mu \begin{cases}
		\trace \left( \bm{\Sigma}_{\text{max}}(\state) \cdot \strain(\dot{\state}) \right) \d \bm{x} & \text{, if \eqref{eqn:ConditionStrainDecomposition}}, \\
			\frac{\varepsilon_{1 2}(\dot{\state}) \cdot \left( \varepsilon_{1 1}(\state)-\varepsilon_{2 2}(\state) \right) - \varepsilon_{1 2}(\state) \cdot \left( \varepsilon_{1 1}(\dot{\state}) - \varepsilon_{2 2}(\dot{\state}) \right)}{\left(\varepsilon_{1 1}(\state)-\varepsilon_{2 2}(\state)\right)^2 + 4 \varepsilon^2_{1 2}(\state)} \\
			\ \ \cdot \trace \left( \bm{\Sigma}_{\text{max}}(\state) \cdot \left( \bm{R}(\state)^\top \strain(\state) \, \bm{Q}(\state) + \bm{Q}(\state)^\top \strain(\state) \, \bm{R}(\state) \right) \right) \\
			\ + \trace \left( \bm{\Sigma}_{\text{max}}(\state) \cdot \bm{Q}(\state)^\top \strain(\dot{\state}) \, \bm{Q}(\state) \right)
		& \text{, else}
	\end{cases} \\
	& - \int_\physdomain \lambda \max(0,\trace(\strain(\state))) \cdot \trace \left( \nabla \state \nabla \vW \right) + \sym\left( \nabla \state \nabla \vW \right) : \stress(\adjstate) \\
	&\hphantom{- \int_\physdomain\,}+ \stress(\state) : \sym\left( \nabla \adjstate \nabla \vW \right) \d \bm{x} \\
	&- \int_\physdomain 2 \mu \begin{cases}
		\trace \left( \bm{\Sigma}_{\text{max}}(\state) \cdot  \sym\left( \nabla \state \nabla \vW \right) \right) \d \bm{x} & \text{, if \eqref{eqn:ConditionStrainDecomposition}}, \\
			\frac{\left( \begin{pmatrix}
					-\varepsilon_{1 2}(\state) & \frac{\varepsilon_{1 1}(\state)-\varepsilon_{2 2}(\state)}{2} \\
					\frac{\varepsilon_{1 1}(\state)-\varepsilon_{2 2}(\state)}{2} & \varepsilon_{1 2}(\state)
				\end{pmatrix} {\nabla \state}^\top \right) : \nabla \vW}{\left(\varepsilon_{1 1}(\state)-\varepsilon_{2 2}(\state)\right)^2 + 4 \varepsilon^2_{1 2}(\state)} \\
			\ \ \cdot \trace \left( \bm{\Sigma}_{\text{max}}(\state) \cdot \left( \bm{R}(\state)^\top \strain(\state) \, \bm{Q}(\state) + \bm{Q}(\state)^\top \strain(\state) \, \bm{R}(\state) \right) \right) \\
			\ +\trace \left( \bm{\Sigma}_{\text{max}}(\state) \cdot \bm{Q}(\state)^\top \sym\left( \nabla \state \nabla \vW \right) \, \bm{Q}(\state) \right) \d \bm{x}
		& \text{, else}
	\end{cases} \\
	&+ \int_\physdomain \Div(\vW) \cdot \left( \frac{\lambda}{2} \max(0,\trace(\strain(\state)))^2 + \mu \trace \left( \bm{\Sigma}_{\text{max}}(\state)^2 \right) + \stress(\state) : \strain(\adjstate) - \nu \right) \d \bm{x} \\
	&+ \int_\shape \frac{G_c}{2} \left(\Div(\vW) - \bm{n}^\top \nabla \vW \bm{n}\right) \d \bm{s}.
\end{align*}
The first integral is zero, as can be seen from the state equation in \eqref{eq:StateEquation}, and the second and third integral together are zero due to the adjoint equation \eqref{eq:AdjointEquation}. The shape derivative is therefore given by \eqref{eqn:ShapeDerivative}.

\bibliographystyle{plainurl} 
\bibliography{literature}

\end{document}